\documentclass{amsart}
\usepackage{graphicx}
\usepackage{amssymb}
\usepackage{epstopdf}
\usepackage{nicefrac}
\usepackage{enumerate}
\usepackage{hyperref}
\usepackage{float}
\usepackage{color}
\usepackage{pst-all}
\usepackage{tabularx}

\usepackage[normalem]{ulem}

\newtheorem{thm}{THEOREM}[section]

\newtheorem{cor}[thm]{COROLLARY}
\newtheorem{defn}[thm]{DEFINITION}
\newtheorem{ex}[thm]{EXAMPLE}

\newtheorem{lemma}[thm]{LEMMA}

\newtheorem{prop}[thm]{PROPOSITION}

\newtheorem{remark}[thm]{REMARK}

\newtheorem{criterion}[thm]{CRITERION}

\newcommand{\ds}{\displaystyle}

\newcommand{\cN}{{\mathcal N}}
\newcommand{\cA}{{\mathcal A}}

\newcommand{\cC}{{\mathcal C}}
\newcommand{\cD}{{\mathcal D}}

\newcommand{\cG}{{\mathcal G}}
\newcommand{\cH}{{\mathcal H}}

\newcommand{\cK}{{\mathcal K}}

\newcommand{\cM}{{\mathcal M}}

\newcommand{\cP}{{\mathcal P}}
 
\newcommand{\cS}{{\mathcal S}}

\newcommand{\cU}{{\mathcal U}}
\newcommand{\cV}{{\mathcal V}}

\newcommand{\diam}{{\rm diam}}

\newcommand{\e}{{\varepsilon}}

\newcommand{\fG}{{\mathfrak{G}}}
\newcommand{\fK}{{\mathfrak{K}}}

\newcommand{\fX}{{\mathfrak{X}}}

\newcommand{\G}{\Gamma}

\newcommand{\Ad}{{\rm Ad}} 
\newcommand{\Aut}{{\rm Aut}} 
\newcommand{\Homeo}{{\rm Homeo}}

\newcommand{\Iso}{{\rm Iso}} 
\newcommand{\mC}{{\mathbb C}}

\newcommand{\mN}{{\mathbb N}}

\newcommand{\mZ}{{\mathbb Z}}

\newcommand{\whg}{\widehat{g}}
\newcommand{\whd}{\widehat{d}}
\newcommand{\whh}{\widehat{h}}

\newcommand{\whC}{\widehat{C}}

\newcommand{\whtau}{{\widehat{\tau}}}

\newcommand{\id}{{\rm id}}

\begin{document}

\title{ Non-Hausdorff germinal groupoids for actions of countable groups}

\author{Olga Lukina}
\address{Olga Lukina, Mathematical Institute, Leiden University, PO Box 9512, 2300 RA Leiden, The Netherlands}
\email{o.lukina@math.leidenuniv.nl}

\thanks{Version date: September 12, 2023. Revision: February 16, 2024}



\begin{abstract}
We study conditions under which the germinal groupoid associated to a minimal equicontinuous action of a countable group on a Cantor set has non-Hausdorff topology. We develop a new criterion, which serves as an obstruction for the \'etale topology on the groupoid to be non-Hausdorff. We use this and other criteria to study the topology of germinal groupoids for a few classes of actions. In particular, we give examples of families of contracting and non-contracting self-similar groups, which are amenable and whose actions have associated germinal groupoids with non-Hausdorff topology.

\end{abstract}

\maketitle
 
\section{Introduction}\label{sec-intro}
 
 In this note, we consider minimal equicontinuous group actions on Cantor sets. We are interested in the question, when is the topology on the associated germinal groupoid   non-Hausdorff?  
 In Section \ref{sec-criteria} we present a list of criteria which allow us to answer this question for a range of specific groupoids.
We then study the Hausdorff property for the topology of actions of specific families of discrete groups. Among these families are iterated monodromy groups of quadratic polynomials, and a family of non-contracting weakly branch groups in \cite{Noce2021}. In these examples, the acting groups are amenable, and many of them have non-Hausdorff associated germinal groupoids.

   An obstruction to the groupoid topology to be non-Hausdorff, for actions of profinite groups,  was developed in \cite{HL2021}. The main technical result of this note,  Theorem~\ref{thm-incl-countable}, extends the known results for profinite group actions to the setting of actions of countable groups. The   proof of Theorem~\ref{thm-incl-countable} uses a variety of techniques developed in the author's previous works on Cantor group actions.

\medskip
 
We now give a detailed description of the problem we consider, and precise statements of our results.
 
Let $\fX$ be a Cantor set, and let $\G$ be an infinite (discrete or profinite) group. Then a \emph{group Cantor action}, or just a \emph{Cantor action}, of $\G$ on $\fX$ is given by a homomorphism $\Phi \colon \G \to \Homeo (\fX)$, and we denote such an action by $(\fX,\G,\Phi)$. Throughout the paper we use the notation $g \cdot x = \Phi(g)(x)$ for the action of $g \in \G$ on $x \in \fX$. We also assume $(\fX,\G,\Phi)$ is minimal, that is, for any $x \in \fX$ the orbit $O(x) = \{g \cdot x \mid g \in \G \}$ is dense in $\fX$.

Let $d$ be a metric on $\fX$ compatible with its topology. The action $(\fX,\G,\Phi)$ is \emph{equicontinuous}, if for any $\epsilon >0$ there exists $\delta>0$ such that for any $g \in \G$ and any $x,y \in \fX$ such that $d(x,y) < \delta$ we have $d(g\cdot x, g \cdot y)< \epsilon$. In this article we are interested in minimal equicontinuous actions on Cantor sets, although some of the theory we discuss holds for a wider class of actions.

Associated to a group action $(\fX,\G,\Phi)$, where $\G$ may be a finite, countable or profinite group, there is a \emph{germinal groupoid} $\cG(\fX,\G,\Phi)$ defined as follows.

\begin{defn}\label{def-groupoid}
Let $(\fX,\G,\Phi)$ be an action of a group $\G$ on a topological space $\fX$.

For $g_1, g_2 \in \G$, we say that $\Phi(g_1)$ and $\Phi(g_2)$ are \emph{germinally equivalent} at $x \in\fX$ if $\Phi(g_1) (x) = \Phi(g_2)(x)$, and there exists an open neighborhood $U \subset \fX$ of $x$ such that the restrictions agree, $\Phi(g_1)|U = \Phi(g_2)|U$. We then write $\Phi(g_1) \sim_x \Phi(g_2)$. 

For $g \in \G$, denote the equivalence class of  $\Phi(g)$ at $x$ by $[g]_x$.
 The collection of germs $\cG(\fX, \G, \Phi) = \{ [g]_x \mid g \in \G ~ , ~ x \in X\}$ with sheaf topology   forms an \emph{\'etale groupoid} modeled on $\fX$. 
 \end{defn}
 
In this paper, we are concerned with the case when $\fX$ is a Cantor set, and $\G$ is a countable group. 
 The motivation to study topological properties of the germinal groupoid $\cG(\fX,\G,\Phi)$ comes from the study of the reduced  $C^*$-algebra $C_r^*(\fX,\G,\Phi)$ associated to a Cantor action  $(\fX,\G,\Phi)$. This $C^*$-algebra is an invariant of the continuous orbit equivalence class of the action, and the study of its $K$-theory offers an approach to the classification of Cantor actions, as used for example in the works \cite{GPS2019,Li2018}. 
 The $C^*$-algebra mentioned above can also be constructed using the  germinal groupoid $\cG(\fX, \G, \Phi)$ associated to the action,    as discussed for example by Renault in \cite{Renault1980}.  In \cite{Renault2008}, Renault assumes that the Cantor action $(\fX,\G,\Phi)$ is topologically free, and thus the germinal groupoid $\cG(\fX, \G, \Phi)$ is a Hausdorff topological space, in order to avoid technical difficulties that arise otherwise. 
 
 If the acting group $\G$ in the minimal equicontinuous action $(\fX,\G,\Phi)$ is abelian, and the action is effective, then it is necessarily topologically free and the associated germinal groupoid is Hausdorff, see \cite[Corollary 2.3]{HL2018}. Thus a germinal groupoid with non-Hausdorff topology can only arise, when the acting group is non-abelian; it is an interesting problem to relate the non-Hausdorff property of the germinal groupoid with the algebraic properties of the acting group. For instance, if $\G$ is nilpotent, or, more generally, has Noetherian property, then the associated germinal groupoid $\cG(\fX,\G,\Phi)$ must have Hausdorff topology, see Section \ref{sec-germ-LQA}, but this question is open for other groups.
 
 An approach to the classification of minimal equicontinuous actions on Cantor sets, using algebraic invariants, was developed in author's works \cite{DHL2017,HL2018,HL2019,HL2021} joint with Dyer and Hurder, see also Section \ref{sec-dirlimitgroups}.  
 In these works, we introduced a notion of complexity of the action, whether they are \emph{stable} or \emph{wild}, which divides the minimal equicontinuous actions of non-abelian groups into classes, see Section \ref{sec-invariants} for details. One can prove that stable actions always have Hausdorff germinal groupoids; for wild actions the situation may vary, for they are divided into further classes of complexity. The least well-behaved class is that of \emph{dynamically wild} actions; we will show that actions with non-Hausdorff germinal groupoid belong to this class. 
The case where $\cG(\fX, \G, \Phi)$  has non-Hausdorff topology may be considered to be  exceptional, and the fact that the topology is non-Hausdorff    has implications for the algebraic structure of   $C_r^*(\fX,\G,\Phi)$, as discussed   in  \cite{BCFS2014,Exel2011}. 

\medskip
We now give the definition of a non-Hausdorff element and discuss the criteria for the germinal groupoid to be non-Hausdorff.

\begin{defn}\label{defn-nonH-rev}
Let $\G$ be a countable or profinite group acting on a topological space $\fX$. Then $g \in \G$ is called a \emph{non-Hausdorff element} at $x \in \fX$ if the following conditions hold: 
\begin{enumerate}
\item $g \cdot x = x$,
\item $g$ is not the identity map on any open neighborhood $W \owns x$, 
\item for any neighborhood $W$ of $x$ there is an open set $O \subset W \subset \fX$ such that $g|O = \id$. 
\end{enumerate}
\end{defn}

Examples of group actions with non-Hausdorff elements include, for instance, the action of the Grigorchuk group in \cite[Section 1.6]{Nekrashevych2005}. Also, the process of the fragmentation of dihedral groups in \cite{Nekrashevych2016} is done by adding non-Hausdorff elements to the group.

In Section \ref{sec-criteria} we show that a germinal groupoid $\cG(\fX,\G,\Phi)$ associated to a group action on a topological space is non-Hausdorff if and only if $\G$ has a non-Hausdorff element as in Definition \ref{defn-nonH-rev}. This is Criterion \ref{crit-nonH-element}, and it follows from the work of Winkelnkemper \cite{Winkel1983}. In some situations, such as Theorems \ref{thm-polyn} and \ref{thm-Md-nonHausd} below, a direct application of this criterion allows to determine that the associated germinal groupoid is non-Hausdorff.

If there is no apparent choice of a non-Hausdorff element in $\G$, then it is usually a very difficult problem to rule out its existence. For this purpose in Section \ref{sec-criteria} we present two criteria which first appeared in \cite{HL2021}. While the first one, Criterion \ref{crit-hausdorff2}, is applicable to actions of discrete or profinite groups, the second one in Theorem \ref{thm-incl-profinite} characterizes the action the profinite group, which is the closure of the action $(\fX,\G,\Phi)$ of a countable group $\G$ in ${\rm Homeo}(\fX)$. The main technical result in this paper, Theorem \ref{thm-main1}, develops an analogue of this criterion for the action of countable groups, obtaining a new invariant of isomorphism classes of minimal equicontinuous actions in the process. 

A question whether a property of a profinite group translates into a property of a dense countable subgroup of this group is very subtle, as their properties may vary widely, especially for non-abelian groups. For instance, it is known that the profinite completion of a nilpotent group is torsion-free if and only if the group itself is torsion free \cite[Corollary 4.7.9]{RZ2000}. On the other hand, Lubotzky \cite{Lubotsky1993} showed that there exist a finitely generated torsion-free linear group whose profinite completion contains torsion elements of any order $r \geq 2$. 
Another related question is the following problem posed by Grothendieck in \cite{Grothen}: let $\G_i$, $i=1,2$ be finitely presented residually finite groups, such that $\G_i$ is dense in a profinite group $\fG_i$, $i=1,2$, and let $h: \G_1 \to \G_2$ be a homomorphism. Let $\widehat h: \fG_1 \to \fG_2$ be the induced homomorphism of profinite groups. If $\widehat h$ is an isomorphism, does this imply that $\G_1$ and $\G_2$ are isomorphic? The first counterexample to this statement in the case when $\G_i$, $i = 1,2$, are finitely generated, was given by Platonov and Tavgen \cite{PT1986}, and later a family of counterexamples was found by Bass and Lubotzky \cite{BL2000}. An uncountable familiy of pairs of finitely generated non-isomorphic groups with isomorphic profinite completions were constructed by Pyber \cite{P2004} and Nekrashevych \cite{Nekr2014}. Finally, finitely presented counterexamples to Grothendieck's conjecture were found by Bridson and Grunewald \cite{BG2004}.

A property of countable group actions, generalizing the notion of a topologically free action, is that such an action is \emph{locally quasi-analytic (LQA)}, see Section \ref{sec-LQA}. For a minimal equicontinuous action $(\fX,\G,\Phi)$ of a countable group $\G$, denote by $\fG(\Phi)$ the closure of the action in ${\rm Homeo}(\fX)$. The closure $\fG(\Phi)$ is a profinite group, and there is an associated action of $\fG(\Phi)$ on $\fX$, denoted by $(\fX,\fG(\Phi),\widehat \Phi)$, see Section \ref{sec-profinite} for details. By the properties of this associated action, minimal equicontinuous actions are divided into \emph{stable} and \emph{wild}, and wild actions are further divided into \emph{wild of finite type}, \emph{wild of flat type}, and \emph{dynamically wild}. The precise definitions of these classes require introduction of a few technical concepts, and we postpone them until Section \ref{sec-invariants}.  By the properties of its closure, a LQA action $(\fX,\G,\Phi)$ of a countable group $\G$ may be stable or wild, and the action which is not LQA is always wild.

Theorem \ref{thm-main1} below proves that a classification similar to the classification of wild actions above holds also for countable groups. The countable analogue of the profinite property of an action being \emph{wild of finite type}, \emph{wild of flat type}, or \emph{dynamically wild}, are the properties that the action is \emph{countably wild of finite type}, \emph{countably wild of flat type}, or \emph{countably dynamically wild}, see Definition \ref{defn-types-wild-countable}. We call the stabilizer $\fG(\Phi)_x$ of the action of the profinite group $\fG(\Phi)$ at $x \in \fX$ the \emph{discriminant group}, and the stabilizer $\cK(\G^x)$ of the action of the countable group $\G$ at $x$ the \emph{kernel} of the action at $x \in \fX$. While the discriminant group is independent of the choice of $x \in \fX$ up to an isomorphism, the kernel $\cK(\G^x)$ depends on the choice of $x \in \fX$. The classification of minimal equicontinuous actions as stable or wild (of different types) is based on the properties of certain subgroups of the discriminant group. A new result in this paper is that, although the kernel $\cK(\G^x)$ depends on $x \in \fX$, certain subgroups of the kernel are preserved. These are the subgroups used to define a new invariant, and develop the classification in Theorem \ref{thm-main1}.

\begin{thm}\label{thm-main1}
Let $(\fX,\G,\Phi)$ be a minimal equicontinuous action of a countable group $\G$ on a Cantor set $\fX$. 
\begin{enumerate}
\item The property that $(\fX,\G,\Phi)$ is countably wild of finite or flat type, or that it is countably dynamically wild is an invariant of the conjugacy class of the action $(\fX,\G,\Phi)$.
\item If $(\fX,\G,\Phi)$ is countably wild of flat type, then the germinal groupoid $\cG(\fX,\G,\Phi)$ has Hausdorff topology (Criterion \ref{crit-upsilon-groups}).
\item If the group $\G$ has a non-Hausdorff element, then $(\fX,\G,\Phi)$ is countably dynamically wild.
\item If $(\fX,\G,\Phi)$ is countably dynamically wild then $(\fX,\G,\Phi)$ is dynamically wild.
\item If $(\fX,\G,\Phi)$ is wild of finite (resp. flat) type, then either $(\fX,\G,\Phi)$ is LQA, or it is countably wild of finite (resp. flat) type.
\end{enumerate}
\end{thm}

Criterion \ref{crit-upsilon-groups} in item (2) of Theorem \ref{thm-main1} does not have a converse, i.e. an action with Hausdorff groupoid can still be countably dynamically wild. Thus Criterion \ref{crit-upsilon-groups} is an obstruction to the existence of non-Hausdorff elements  in $\G$;  if this criterion is not satisfied, then the groupoid $\cG(\fX,\G,\Phi)$ may or may not have non-Hausdorff topology. We present a counterexample in Theorem \ref{thm-no-crit} below, referring the reader to Section \ref{sec-trees} for necessary background on actions on rooted trees.

\begin{thm}\label{thm-no-crit}
There exists a contracting self-similar group $\G$ acting on the boundary of a binary tree $T$, such that the action $(\partial T, \G,\Phi)$ does not satisfy Criterion \ref{crit-upsilon-groups} while its germinal groupoid $\cG(\partial T, \G,\Phi)$ has Hausdorff topology.
\end{thm}

We present two families of examples of group actions on rooted trees for which the associated germinal groupoids have Hausdorff or non-Hausdorff topology. The first family is a class of contracting self-similar actions, for which there is a specific Criterion \ref{crit-contracting}, which allows to rule out the existence of non-Hausdorff elements using torsion. This criterion appeared earlier in author's work \cite{Lukina2021}. Contracting self-similar groups arise, for instance, as iterated monodromy groups associated to complex polynomials, and we have the following result.

Let $f(x)$ be a polynomial of degree $d = 2$ over $\mathbb{C}$. Denote by $c$ the critical point of $f(x)$, and by $P_c = \{f^m(c) \mid m \geq 1\}$ the orbit of the critical point. Then $f(x)$ is \emph{post-critically finite (PCF)}, if $P_c$ is finite. If $P_c$ is finite, then either it consists of a single periodic cycle, which includes the critical point $c$, or $P_c$ consist of a single strictly pre-periodic orbit, i.e. there exist integers $k,m \geq 1$ such that $f^{n+k}(x) = f^n(x)$ for all $n \geq m$.

Associated to $f(x)$, there is a discrete group ${\rm IMG}(f)$ called the \emph{discrete iterated monodromy group} of $f(x)$, which acts on preimages of a non-critical point $t \in \mC$ by permutations. Such preimages are vertices of a binary tree $T$ with root $t$, see Section \ref{sec-trees} for the background on actions on rooted trees.

\begin{thm}\label{thm-polyn}
Let $f(x)$ be a quadratic PCF polynomial, and let $\cG(\partial T, {\rm IMG}(f),\Phi)$ be the germinal groupoid associated to the action of the discrete iterated monodromy group ${\rm IMG}(f)$ on the boundary $\partial T$ of a binary tree $T$. Then the germinal groupoid $\cG(\partial T, {\rm IMG}(f),\Phi)$ has non-Hausdorff topology if and only if the post-critical orbit $P_c$ is strictly pre-periodic, and $P_c$ has cardinality at least $3$. Otherwise $\cG(\partial T, {\rm IMG}(f),\Phi)$ has Hausdorff topology.
\end{thm}
 
Theorem \ref{thm-polyn} is a consequence of a more general result proved in Theorems  \ref{thm-hausdorff-periodic} and \ref{thm-thm12}, for groups $\fK(v)$ and $\fK(w,v)$ defined in Sections \ref{sec-kv} and \ref{sec-kwv} respectively. Discrete iterated monodromy groups of quadratic polynomials form a proper subset of the set of these groups, see \cite{BN2008} for details. A similar result for a single specific choice of $v$ and $wv$ was obtained in \cite{Lukina2021}. Theorem \ref{thm-polyn} does not follow from the results of \cite{Lukina2021}, since, as discussed above, the properties of profinite completions and closures of actions need not translate into propertes of discrete groups. The proof of Theorem \ref{thm-hausdorff-periodic} in Section \ref{sec-kv} is by a direct computation.

Criterion \ref{crit-contracting} only applies to actions of self-similar groups which are contracting. Our next result shows that the contracting property is crucial here, as there exist self-similar weakly branch groups which are not contracting, and whose germinal groupoids are non-Hausdorff with non-Hausdorff elements of infinite order.

Noce \cite{Noce2021} constructed a family of finitely generated weakly branch non-contracting groups $\cM(d)$, where $\cM(d)$ acts on a $d$-ary tree, for every $d \geq 2$. We describe this family of actions in Section \ref{sec-Md}, where we also prove the following result.

\begin{thm}\label{thm-Md-nonHausd}
For any $d \geq 3$, the germinal groupoid $\cG(\partial T, \cM(d),\Phi)$, associated to the action of a weakly branch non-contracting group $\cM(d)$ of automorphisms of a $d$-ary tree $T$, has non-Hausdorff topology. 
\end{thm}

Groups in Theorem \ref{thm-polyn} are generated by bounded automata, and groups in Theorem \ref{thm-Md-nonHausd} are generated by automata with linear activity, and so both families consist of amenable groups, see Remark \ref{remark-bounded-amenable} and Theorem \ref{exp-activity} for details. Thus, along with the Grigorchuk group which is amenable and has non-Hausdorff germinal groupoid, Theorems \ref{thm-polyn} and \ref{thm-Md-nonHausd} provide examples of classes of amenable groups with non-Hausdorff germinal groupoid, see also Remark \ref{remark-amenable-nonHausd}.


The rest of the paper is organized as follows. In Section \ref{sec-invariants} we recall the basic properties of minimal equicontinuous actions, and define their classifying algebraic invariants. In Section \ref{sec-criteria} we discuss Criteria \ref{crit-nonH-element} and \ref{crit-hausdorff2}, and prove Theorem \ref{thm-main1}, which gives us Criterion \ref{crit-upsilon-groups}. In Section \ref{sec-trees} we recall the necessary background on actions of groups of automorphisms of rooted trees, and state Criterion \ref{crit-contracting}. In Section \ref{sec-PCF} we prove Theorems  \ref{thm-no-crit} and \ref{thm-polyn}, and in Section \ref{sec-Md} we prove Theorem \ref{thm-Md-nonHausd}.

\section{Algebraic invariants for equicontinuous actions}\label{sec-invariants}

In this section we give a brief outline of theory of algebraic classifying invariants for minimal equicontinuous group actions, developed in the author's works \cite{DHL2017,HL2019,HL2021} joint with Dyer and Hurder. This is needed in order to introduce and prove the criteria for the non-Hausdorff property in Section \ref{sec-criteria}.

\subsection{Adapted sets and neighborhood bases}\label{sec-adapted}

Let $\fX$ be a Cantor set, that is, a compact totally disconnected metrizable space without isolated points. Let $\G$ be a countable group, and
let $(\fX,\G,\Phi)$ be a minimal equicontinuous Cantor action, as defined in Section \ref{sec-intro}.

Let ${\rm CO}(\fX)$ denote the collection  of all clopen (closed and open) subsets of  $\fX$, which forms a basis for the topology of $\fX$. 
For any $\phi \in \Homeo(\fX)$ and   any $U \in {\rm CO}(\fX)$, the image $\phi(U) \in {\rm CO}(\fX)$.  
The following   result is folklore, and a proof is given in \cite[Proposition~3.1]{HL2018}.
 \begin{prop}\label{prop-CO}
 For $\fX$ a Cantor set, a minimal   action   $\Phi \colon \G \times \fX \to \fX$  is  equicontinuous  if and only if  the $\G$-orbit of every $U \in {\rm CO}(\fX)$ is finite for the induced action $\Phi_* \colon \G \times {\rm CO}(\fX) \to {\rm CO}(\fX)$.
\end{prop}

Proposition \ref{prop-CO} leads to the existence of certain dynamically defined subsets, which we call \emph{adapted}. The techniques we use to study minimal equicontinuous Cantor actions rely on the existence of such adapted sets.

\begin{defn}\label{defn-adaptedset}
    We say that a subset $U \subset \fX$ is \emph{adapted} to a Cantor action $(\fX,\G,\Phi)$ if for every $g \in \G$ we have either $g \cdot U = U$, or $g \cdot U \cap U = \emptyset$. 
\end{defn}

The proof of \cite[Proposition~3.1]{HL2018} implies that, given a clopen $W \subset \fX$ and $x \in W$,  there is an adapted clopen set $U$ with $x \in U \subset W$. Then  the set of ``return times'' to $U$, 
 \begin{equation}\label{eq-adapted}
\G_U = \left\{g \in \G \mid g \cdot U  \cap U \ne \emptyset  \right\} ,
\end{equation}
is a subgroup of   $\G$, called the \emph{stabilizer} of $U$, or the \emph{isotropy subgroup} of $\G$ at $U$.   The  translates $\{ g \cdot U \mid g \in \G\}$ form a finite clopen partition of $\fX$, and are in 1-1 correspondence with the quotient space $X_U = \G/\G_U$. The group $\G$ acts by permutations of the finite set $X_U$, and so the stabilizer group $\G_U \subset \G$ has finite index.  If $V \subset U$ is a proper inclusion of adapted sets, then the inclusion $\G_V \subset \G_U$ is also proper. Thus Proposition \ref{prop-CO} implies the following result.

\begin{prop}\label{prop-adaptedchain}
Let  $(\fX,\G,\Phi)$   be a minimal equicontinuous Cantor    action. Given $x \in \fX$, there exists a properly descending chain of clopen sets $\cU = \{U_{\ell} \subset \fX  \mid \ell \geq 0\}$, $U_0 = \fX$, such that
    $x \in U_{\ell +1} \subset U_{\ell}$  is a proper inclusion for all $ \ell \geq 0$,   $\cap_{\ell \geq 0}  \ U_{\ell} = \{x\}$, and  each $U_{\ell}$ is adapted to the action $\Phi$. 
 \end{prop}

A descending chain of adapted subsets, as in Proposition \ref{prop-adaptedchain}, is called an \emph{adapted neighborhood basis} at $x \in \fX$ for the action $\Phi$. 

By the remarks before Proposition \ref{prop-adaptedchain}, associated to an adapted neighborhood basis $\cU$ at $x \in \fX$ there is a chain 
 $$\G^x_\cU: \G = \G_0 \supset \G_1 \supset \G_2 \supset \cdots$$ of finite index subgroups of $\G$. Denote by $X_\ell = \G/\G_\ell$ the quotient space, and by $p^{\ell+1}_\ell: X_{\ell+1} \to X_\ell: g \G_{\ell+1} \mapsto g\G_\ell$ the coset inclusions. 
  The correspondence of the quotient space $X_\ell$ with the set of translates of $U_\ell$ defines a projection $f_\ell: \fX \to X_\ell$, where $f_\ell(x) = g \G_\ell$ if and only if $x \in g \cdot U_\ell$. The projections are equivariant with respect to the action of $\G$ on $\fX$ and the left action of $\G$ on the coset space $X_\ell$, that is, we have $g \cdot f_\ell(x) = f_\ell (g \cdot x)$. Since $\cap_{\ell \geq 0} U_\ell = \{x\}$, there is a homeomorphism
  \begin{align}\label{eq-invlimit}f_\infty: \fX \to X_\infty = \lim_{\longleftarrow}\{p^{\ell+1}_\ell : X_{\ell+1} \to X_\ell \mid \ell \geq 0\} = \{(f_0(x),f_1(x),\ldots)\} \subset \prod_{\ell \geq 0} X_\ell,\end{align}
equivariant with respect to the given action of $\G$ on $\fX$, and the induced action of $\G$ on $X_\infty$ defined by
  \begin{align}\label{eq-invlimit-action}\G \times X_\infty \to X_\infty: (g, (x_0,x_1,\ldots)) \mapsto (g \cdot x_0, g \cdot x_1,\ldots).\end{align}
The equations \eqref{eq-invlimit}-\eqref{eq-invlimit-action} provide an \emph{inverse limit representation} for a minimal equicontinuous action $(\fX,\G,\Phi)$. Such a representation is not unique and depends on the choice of a point $x \in \fX$ and of an adapted neighborhood basis $\cU$. However, all such representations are conjugate, see \cite{DHL2016a} for details.

\subsection{Locally quasi-analytic actions} \label{sec-LQA}

The notion of a \emph{quasi-analytic} topological action of a   pseudogroup  on a connected topological space was introduced by Haefliger \cite{Haefliger1985}. 
  {\'A}lvarez L{\'o}pez and Candel in   \cite[Definition~9.4]{ALC2009}, and  later 
   {\'A}lvarez L{\'o}pez and Moreira Galicia in \cite[Definition~2.18]{ALM2016}, adapted the notion of a quasi-analytic topological action to the more general case where the action space need not be connected.

  \begin{defn} \cite[Definition~9.4]{ALC2009} \label{def-LQA} A Cantor action       $(\fX,\G,\Phi)$  is   \emph{locally quasi-analytic}, or simply   \emph{LQA}, if there exists $\e > 0$ such that for any adapted   set $U \subset \fX$ with $\diam (U) < \e$,  and  for any adapted subset $V \subset U $, and elements $g_1 , g_2 \in \G$
 \begin{equation}\label{eq-lqa}
  \text{if the restrictions} ~~ \Phi(g_1)|V = \Phi(g_2)|V, ~ \text{ then}~~ \Phi(g_1)|U = \Phi(g_2)|U.   
\end{equation}
If  \eqref{eq-lqa} holds for $U=\fX$, then the action of $\G$ is \emph{quasi-analytic}.
\end{defn}

 In other words, an action $(\fX,\G,\Phi)$ is locally quasi-analytic on open sets of diameter $\varepsilon >0$, if for every $g \in \G$ the homeomorphism $\Phi(g) $ extends uniquely from small open sets to sets of diameter at least $\varepsilon$. Thus the property of local quasi-analyticity quantifies the complexity of an action.

The notion of a locally quasi-analytic action can be seen as a generalization of the notion of a topologically free action, as we explain now.
Recall that an action $(\fX,\G,\Phi)$ is \emph{effective} if $g \cdot x = x$ for all $x \in \fX$ implies $g = \id$.

\begin{defn}\label{defn-topfree}
A Cantor action $(\fX,\G,\Phi)$ is said to be \emph{topologically free}   if the set 
  $$\Iso(\fX,\G,\Phi) = \{ x \in \fX \mid g\cdot x = x \textrm{ for some }g \in \G\} $$ is meager   in $\fX$, 
  and an action $(\fX,\G,\Phi)$ is said to be \emph{free} if $\Iso(\fX,\G,\Phi) $ is an empty set.
  \end{defn}
  
Note that if $e \ne g \in \G$ and $\Phi(g)$ acts trivially on $\fX$, then $\Iso(\fX,\G,\Phi)  = \fX$, and thus a topologically free action must be effective.  The relation between quasi-analytic and topologically free actions is given by the following proposition.

\begin{prop}\cite[Proposition 2.2]{HL2018}
An effective group Cantor action $(\fX,\G,\Phi)$ is quasi-analytic if and only if it is topologically free.
\end{prop}

If the group $\G$ is abelian, it is an exercise to show that an effective minimal Cantor action $(\fX,\G,\Phi)$ must be topologically free; see for instance \cite[Corollary~2.3]{HL2018}. Examples of equicontinuous Cantor actions which are locally quasi-analytic, but not quasi-analytic,  are easily constructed, see \cite{DHL2017,HL2019}.
 
To clarify the relation between a topologically free and a locally quasi-analytic minimal equicontinuous actions, let $U$ be an adapted set for $(\fX, \G,\Phi)$, and consider the restricted homomorphism 
 $$\Phi|_U: \G_U \to {\rm Homeo}(U),$$ 
 where $\G_U$ is the isotropy group at $U$ defined by \eqref{eq-adapted}. 
 If $(\fX,\G,\Phi)$ is not topologically free, then $\Phi|_U$ may have a non-trivial kernel. If $(\fX, \G,\Phi)$ is locally quasi-analytic on sets of diameter $\varepsilon >0$, and $\diam (U) < \varepsilon$, then the induced action 
  $$\widetilde{\Phi}|_U: \G_U/\ker \Phi_U \to {\rm Homeo}(U)$$ 
of $\G_U$ on $U$  is topologically free.

\subsection{Profinite closure of the action} \label{sec-profinite}

In order to define the direct limit group invariants of a minimal equicontinuous action $(\fX,\G,\Phi)$, we introduce a group chain model for the closure of the action. Namely, let
$$\fG(\Phi) = \overline{\Phi(G)} \subset {\rm Homeo}(\fX),$$
be the closure of the action in the uniform topology. The closure is a profinite group. Elements of $\fG(\Phi)$ are the sequences $\widehat{g} = \{g_\ell \mid g_\ell \in \G, \ell \geq 0\}$, and the image  $\Phi(\G)$ is a dense subgroup of $\fG(\Phi)$. If the action $(\fX,\G,\Phi)$ is effective, $\Phi(\G)$ is identified with $\G$. In any case, the action $\Phi$    induces an action $\widehat{\Phi}:\fG(\Phi) \to {\rm Homeo}(\fX)$, which is transitive since $(\fX,\G,\Phi)$ is assumed to be minimal. 

\begin{remark}
{\rm
In topological dynamics $\fG(\Phi)$ is also known as the \emph{Ellis (semi-)group} of the action $(\fX,\G,\Phi)$, see \cite{Auslander1988}. The Ellis semigroup of a (not necessarily equicontinuous) group action $(\fX,\G,\Phi)$ is defined as the closure of the action in the topology of pointwise convergence. In our case, since $(\fX,\G,\Phi)$ is equicontinuous, this topology coincides with the uniform topology, and the Ellis semi-group is a group. Thus the closure $\fG(\Phi)$ can also be referred to as the Ellis group of the action $(\fX,\G,\Phi)$. 
}
\end{remark}

For a point $x \in \fX$, denote by 
  \begin{align*}\fG(\Phi)_x = \{\widehat{g} \in \fG(\Phi) \mid \widehat{g} \cdot x =x\}\end{align*}
the stabilizer, or the isotropy group of the action of $\fG(\Phi)$ at $x \in \fX$, called the \emph{discriminant group} of $(\fX,\G,\Phi)$.  The quotient $\fG(\Phi)/\fG(\Phi)_x$ is homeomorphic to $\fX$, and is a homogeneous space for the action of $\fG(\Phi)$. 
The action of $\G$ on the coset space $\fX \cong \fG(\Phi)/\fG(\Phi)_x$ gives a \emph{homogeneous model} for the minimal equicontinuous action $(\fX,\G,\Phi)$.

\begin{defn}\label{defn-wild}
A minimal equicontinuous action $(\fX,\G,\Phi)$ is \emph{stable} if the action of its closure $\fG(\Phi) = \overline{\Phi(\G)} \subset {\rm Homeo}(\fX)$ on $\fX$ is locally quasi-analytic, and otherwise $(\fX,\G,\Phi)$ is wild.
\end{defn}

Given a group chain $\G^x_{\cU}$ for an adapted neighborhood system $\cU$, for each $\ell \geq 0$, let $C_\ell = \bigcap_{g \in \G} g \G_\ell g^{-1}$ be the normal core of $\G_\ell$. Then $\{C_\ell \mid \ell \geq 0\}$ is a descending chain of finite index normal subgroups of $\G$. The proof of the following theorem can be found, for instance, in \cite[Theorem 4.4]{DHL2016a}.

\begin{thm}\label{thm-quotientspace}
Let    $(\fX,\G,\Phi)$ be a Cantor action, and suppose that $\G^x_\cU = \{\G_{\ell}\mid \ell \geq 0 \}$ is the group chain associated to an adapted neighborhood basis $\cU$ at $x \in \fX$.  Let $\{C_\ell \mid \ell \geq 0\}$ be the chain of normal cores of the groups in $\G^x_\cU$. Then there is an isomorphism of topological groups
\begin{align}\label{eq-Cspace}\widehat f \colon    \fG(\Phi) \to \whC_{\infty} = \lim_{\longleftarrow}\{\G/C_{\ell+1} \to \G/C_\ell \mid \ell \geq 0\},\end{align}  
where $\G/C_{\ell+1} \to \G/C_\ell$ are the maps on cosets induced by the inclusions $C_{\ell+1} \to C_\ell$. Moreover, \eqref{eq-Cspace} restricts to the isomorphism of the  discriminant groups
  $$\widehat f \colon \fG(\Phi)_x \to \cD_x = \lim_{\longleftarrow}\{\G_{\ell+1}/C_{\ell+1} \to \G_\ell/C_\ell \mid \ell \geq 0\}.$$
\end{thm}

Elements of $\whC_\infty$ are chains of cosets $\whg = \{g_\ell C_\ell \mid \ell \geq 0\}$ such that $g_{\ell+1}C_{\ell+1} \subset g_\ell C_\ell$. In the rest of the paper we use the notation $\whg = (g_\ell)$ whenever this does not lead to confusion.

\subsection{Centralizer and stabilizer  direct limit groups}\label{sec-dirlimitgroups}

By Theorem \ref{thm-quotientspace} we have $X_\infty \cong \whC_\infty / \cD_x$.

We will now define direct limit stabilizer and centralizer groups associated to a minimal equicontinuous action $(\fX,\G,\Phi)$. We assume that $\G$ acts effectively on $\fX$, so that $\G$ injects onto a dense subgroup of $\fG(\Phi)$, also denoted by $\G$. Then the image $\widehat f(\G) \subset \whC_\infty$, also denoted by $\G$, is dense in $\whC_\infty$. For each $\ell \geq 0$, the topological closure  $\widehat \G_\ell  = \overline{\G_\ell}$ in $\whC_\infty$ is a clopen subgroup, in particular,  it has finite index in $\whC_\infty$. From the homogeneous model, $U_\ell = \widehat \G_\ell /\cD_x$.

The action of the discriminant group $\cD_x$ on the homogeneous space $X_\infty$ can be seen as the adjoint action, namely, for $\widehat h \in \cD_x$ and $\widehat g \in \whC_\infty$
  \begin{align}\label{eq-adjoint}\Ad(\widehat h) (\widehat g \, \cD_x) = \widehat h \, \widehat g \, \widehat h^{-1} \cD_x =  \widehat h \,\widehat g \,  \cD_x.\end{align}
 The adjoint action of $\whh$ on $\whg \, \cD_x$ is trivial, i.e. 
   $$\Ad(\widehat h)(\widehat g \, \cD_x) = \widehat g \, \cD_x,$$ 
if and only if the commutator $[\widehat h^{-1} ,\widehat g^{-1}]$ is in $\cD_x$. By assumption $\widehat h \in \cD_x$, so $\whh$ fixes $x \in U_\ell$. Then the action of $\whh$ preserves $U_\ell$, since $U_\ell$ is adapted, and for $\widehat g \in \widehat \G_\ell$ we have $\Ad(\widehat h)(\widehat g) \in \widehat \G_\ell$, i.e. the adjoint action of $\widehat h$ preserves the clopen subgroups $\widehat \G_\ell$, for $\ell \geq 0$.

To define the direct limit group invariants, we distinguish between the elements of $\cD_x$ that commute with all elements in $\widehat \G_\ell$, and those whose commutators with elements of $\widehat \G_\ell$ are in $\cD_x$ but may be non-trivial. By the above remarks in both cases the action on  $U_\ell$ via \eqref{eq-adjoint} is trivial.

We define an increasing chain of \emph{stabilizer subgroups} of $\cD_x$ by
  \begin{align}\label{eq-Kchain} K(\Phi) = \{K_\ell\}_{ \ell \geq 0}, &  &K_\ell =  \left\{ \widehat{h} \in \cD_x \mid  [\widehat{h}, \widehat{g}] \, \in \cD_x, \, \textrm{ for all }\widehat{g} \in \widehat{\G}_\ell \right\},\end{align}
that is, $K_\ell$ consists of all elements in $\cD_x$ which act  on the clopen set $U_\ell$ by the identity map. Indeed, we have the homogeneous model $U_\ell = \widehat \G_\ell / \cD_x$, i.e. every point in $x \in U_\ell$ corresponds to a coset $\widehat g \, \cD_x$ for some $\widehat g \in \widehat \G_\ell$. Suppose $\widehat h \in K_\ell$ as in \eqref{eq-Kchain}. Since $K_\ell$ and $\widehat \G_\ell$ are groups, then for any $\widehat g \in \widehat \G_\ell$ we have $[\widehat h^{-1}, \widehat g^{-1}] \in \cD_x$, that is, 
  $$\widehat h^{-1} \widehat g^{-1} \widehat h \, \widehat g \, \cD_x = \cD_x,$$
and, multiplying by $\widehat g\, \widehat h$ on the left we obtain $\widehat h \, \widehat g \, \cD_x = \widehat g \, \widehat h \, \cD_x = \widehat g \, \cD_x$. Since this holds for any $x = \widehat g\, \widehat \G_\ell$, we obtain that $\widehat h$ acts as the identity on $U_\ell$. Conversely, if $\widehat h \in \cD_x$ acts as the identity on $U_\ell$, then for any $\widehat g^{-1} \in \widehat \G_\ell$ we have 
   $$\widehat h \, \widehat g^{-1} \cD_x = \widehat g^{-1}\cD_x,$$
and so $\widehat h^{-1} \, \cD_x = \cD_x = \widehat g \, \widehat h^{-1} \, \widehat g^{-1} \cD_x$, which implies $[\widehat h, \widehat g] \in \cD_x$.
We note that if an action $(\fX,\G,\Phi)$ is wild, then $K_\ell$ is non-trivial for all $\ell \geq 0$.

Define another increasing chain of \emph{centralizer subgroups} of $\cD_x$ by 
 \begin{align}\label{eq-Zchain} Z(\Phi) = \{Z_\ell\}_{ \ell \geq 0}, & & Z_\ell = \left\{ \widehat{h} \in \cD_x \mid [\widehat{h},\widehat{g}] = {\id} \in \cD_x \, \textrm{ for all }\widehat{g} \in \widehat{\G}_\ell \right\},\end{align}
that is, $Z_\ell$ contains the elements in $\cD_x$ whose adjoint action on the clopen subgroup $\widehat \G_\ell$ is trivial, i.e. the elements of $Z_\ell$ are in the centralizer of $\widehat \G_\ell$. 

Clearly we have the inclusions $Z_\ell \to K_\ell$. Also, since $U_{\ell +1} \subset U_\ell$, there are inclusions $\iota_\ell^{\ell+1}: K_\ell \to K_{\ell+1}$ for $\ell \geq 0$, which restrict to the inclusions $\iota_\ell^{\ell+1}: Z_\ell \to Z_{\ell+1}$. Taking the direct limits of the group chains $K(\Phi)$ and $Z(\Phi)$ with respect to the inclusions $\iota^{\ell+1}_\ell$, we obtain the \emph{ (direct limit) stabilizer} and the \emph{centralizer groups}
  \begin{align}\label{eq-directlimk} \Upsilon^x_s(\Phi) =\lim_{\longrightarrow} 
\cS(K_\ell, \iota_\ell^{\ell+1},\mN), \quad \textrm{ and } \quad \Upsilon^x_c(\Phi) =\lim_{\longrightarrow} 
\cS(Z_\ell, \iota_\ell^{\ell+1},\mN),
 \end{align}
where indices $\ell$ run through the natural numbers $\mN$.  

 The limit groups \eqref{eq-directlimk} were defined in \cite{HL2021}, where it was proved that, although the groups in \eqref{eq-directlimk}  may depend on the choice of the adapted neighborhood basis $\cU$ (which determines the group chain $\G^x_\cU$, and therefore the system of neighborhoods of the identity $\{\widehat{C}_\ell\}_{\ell \geq 0}$), the isomorphism classes of the direct limits do not depend on these choices, and so we have:

\begin{thm}[{\cite[Theorem 4.15]{HL2021}}]\label{thm-isomupsilon}
Let $(\fX,\G,\Phi)$ be a minimal equicontinuous action of a finitely generated group $\G$ on a Cantor set $\fX$. Then the direct limit isomorphism classes $\Upsilon_s(\Phi)$ and $\Upsilon_c(\Phi)$ of the groups $\Upsilon^x_s(\Phi)$ and $\Upsilon^x_c(\Phi)$ are invariants of the conjugacy class of the action $(\fX,\G,\Phi)$.
\end{thm}

The groups $\Upsilon_c(\Phi)$ and $\Upsilon_s(\Phi)$ are used in the classification of minimal equicontinuous Cantor actions as we describe below, see \cite{HL2021} for details. Before we outline this theory, we recall a few notions which allow us to compare the direct limit groups.

\begin{defn}
The direct limit group $ {\displaystyle \Upsilon^x_s(\Phi) =\lim_{\longrightarrow} 
\cS(K_\ell, \iota_\ell^{\ell+1},\mN)}$ (and similarly for $\Upsilon_c^x(\Phi)$) is \emph{bounded}, if there is $k \in \mN$ such that for all $\ell' \geq \ell \geq k$, the maps $\iota^{\ell'}_\ell = \iota^{\ell'}_{\ell' - 1} \circ \cdots \iota^{\ell+1}_\ell$ are group isomorphisms.
\end{defn}

Recall from Definition \ref{defn-wild} that $(\fX,\G,\Phi)$ is stable if and only if the action of the profinite group $\fG(\Phi)$ on $\fX$ is locally quasi-analytic. It was proved in \cite[Theorem 5.3]{HL2021} that the action $(\fX,\G,\Phi)$ is stable if and only if the stabilizer group $\Upsilon_s(\Phi)$ is bounded.

If the action of $\fG(\Phi)$ on $\fX$ is locally quasi-analytic, then the action of the countable dense subgroup $\G \subset \fG(\Phi)$ on $\fX$ is also locally quasi-analytic. The converse need not hold. Indeed, \cite[Theorem 1.4]{HL2023} constructs an uncountable family of pairwise non-conjugate topologically free actions of the Heisenberg group $\cH$, such that the action of the profinite group $\overline{\Phi(\cH)}$ is not locally quasi-analytic.

The proofs of the following standard facts can be found, for instance, \cite[Chapter VIII, Section 2]{EilenbergSteenrod1952}.

 \begin{defn}\label{def-directedgroupmorphism}
  A map $\Xi$ between directed systems of groups 
  $\cS(G_\ell, \phi^{\ell'}_{\ell} , \mN)$  
  and $\cS(H_{k}, \psi^{k'}_{k} , \mN)$
 is an order-preserving map $\xi \colon \mN \to \mN$, and for each $\ell \in \mN$ 
 a group homomorphism $\xi_{\ell} \colon G_{\ell} \to H_{\xi(\ell)}$ such that for $\ell < \ell'$ we have
 $$   \xi_{\ell'} \circ \phi^{\ell'}_{\ell} = \phi^{\xi(\ell')}_{\xi(\ell)} \circ \xi_{\ell} \colon G_{\ell} \to H_{\xi(\ell')} ~ .$$
\end{defn}

\begin{prop}\label{prop-directedgroupisomorphism} 
 A map $\Xi$ between directed systems of groups 
  $\cS(G_{\ell}, \phi^{\ell'}_{\ell} , \mN)$  
  and $\cS(H_{k}, \psi^{k'}_{k} , \mN)$ induces a homomorphism 
$\ds \underrightarrow{\Xi} \colon \varinjlim \cS(G_{\ell}, \phi^{\ell'}_{\ell} , \ell) \longrightarrow \varinjlim \cS(H_{k}, \psi^{k'}_{k} , \mN)$ of the direct limit groups. 

If each $\xi_{\ell} \colon G_{\ell} \to H_{\xi(\ell)}$ for $\ell \in \mN$ is a  monomorphism of groups, then the induced map $\ds  \underrightarrow{\Xi}$ of the direct limit groups is a group monomorphism.

If each $\xi_{\ell} \colon G_{\ell} \to H_{\xi(\ell)}$ for $\ell \in \mN$ is a  isomorphism of groups, then the induced map $\ds  \underrightarrow{\Xi}$ of the direct limit groups in a group isomorphism.
\end{prop}

A subset $\Lambda \subset \mN$ of a directed set is said to be \emph{cofinal} if for each $\ell \in \mN$, there exists $\ell' \in \Lambda$ with $\ell < \ell'$. Then we have:
\begin{prop}\label{prop-directedgroupcofinal}
Let  $\cS(G_{\ell}, \phi^{\ell'}_{\ell} , \mN)$  be a directed systems of groups, and $\Lambda \subset \mN$ be a cofinal set. Then the inclusion $\Lambda \to \mN$ induces a group isomorphism
$\ds  \varinjlim \cS(G_{\ell}, \phi^{\ell'}_{\ell} , \Lambda) \cong   \varinjlim \cS(G_{\ell}, \phi^{\ell'}_{\ell} , \mN)$.
\end{prop}

The properties of the groups $\Upsilon_c(\Phi)$ and $\Upsilon_s(\Phi)$ allow us to develop further classification of wild minimal equicontinuous actions as in the following definition, see \cite{HL2021} for details.

\begin{defn}\label{defn-types-wild}
Let   $(\fX,\G,\Phi)$ be a minimal equicontinuous Cantor action which is wild. Then $(\fX,\G,\Phi)$ is:
\begin{enumerate}
\item \emph{wild of finite type}  if    the stabilizer  group $\Upsilon_s(\Phi)$   is unbounded, and represented by a chain of finite groups;
\item \emph{wild of flat type}  if    the stabilizer  group $\Upsilon_s(\Phi)$   is unbounded, and $\Upsilon_c(\Phi) = \Upsilon_s(\Phi)$; 
\item \emph{dynamically wild}  if    the stabilizer  group $\Upsilon_s(\Phi)$   is unbounded, and the action is not of flat type.
\end{enumerate}
\end{defn} 

\begin{thm}\cite{HL2021}\label{HL2021}
All possibilities in Definition \ref{defn-types-wild} are realized, that is, there exist minimal equicontinuous actions which are wild of finite or flat type, and which are dynamically wild.

\end{thm}

Further examples of wild Cantor actions with various values of the direct limit centralizer and stabilizer groups can be found in the author's work \cite{ALBLLN2020}, joint with \'Alvarez L\'opez, Barral L\' ijo and Nozawa.  

\begin{ex}\label{ex-trivialZ}
{\rm
Let $(\fX,\G,\Phi)$ be a minimal equicontinuous action, and suppose the profinite group $\fG(\Phi)$ is isomorphic to the wreath product of finite subgroups $\{G_\ell \subset Sym(m_\ell) \mid  \ell \geq 1\}$, where $Sym(m_\ell)$ denotes the symmetric group, and such that $G_\ell$ acts transitively on the set of $m_\ell$ symbols, for $\ell \geq 1$. By \cite[Theorem 1.6]{ALBLLN2020} such an action is dynamically wild, with non-trivial profinite stabilizer group $\Upsilon_s(\Phi)$, and the trivial profinite centralizer group $\Upsilon_c(\Phi)$. }
\end{ex}

\section{Criteria for the non-Hausdorff property}\label{sec-criteria}

In this section we compose a list of criteria which allow to detect whether the germinal groupoid associated to a given minimal equicontinuous action is Hausdorff. Criteria \ref{crit-nonH-element} and \ref{crit-hausdorff2} previously appeared in \cite{HL2021} and they apply to actions of countable or profinite groups. Criterion \ref{crit-upsilon-groups} is new. It is an obstruction to the existence of a non-Hausdorff element for actions of countable groups, which is more sensitive than a similar criterion for profinite groups in \cite{HL2021}.

\subsection{Germinal groupoid and non-Hausdorff elements}\label{subsec-nonHaus}  

The germinal groupoid $\cG(\fX,\G,\Phi)$ associated to an action of a group $\G$ on a topological space $\fX$ was defined in Definition \ref{def-groupoid}. A non-Hausdorff element $g \in \G$ was defined in Definition \ref{defn-nonH-rev}. In this section, $\fX$ is any topological space, and $\G$ is a countable or profinite group acting on $\fX$.

We recall the following result.
 \begin{prop}\cite[Proposition 2.1]{Winkel1983}\label{prop-hausdorff}
The germinal groupoid $\cG(\fX, \G, \Phi)$ is Hausdorff at $[g]_x$  if and only if, for all $[g']_x \in \cG(\fX, \G, \Phi)$ with $g \cdot x = g' \cdot x = y$, if there exists a sequence   $\{x_n\} \subset \fX$ which converges to $x$, and such that $[g]_{x_n} = [g']_{x_n}$ for all $n$, then $[g]_{x} = [g']_{x}$. 
\end{prop}

In Proposition \ref{prop-hausdorff}, consider the composition of maps $h = g^{-1}\circ g'$. Since $g \cdot x = g' \cdot x$, then $h \cdot x = g^{-1}  g'\cdot x = x$. Denote by $[\id]_x$ the germ of the identity map at $x \in \fX$. Then the statement of Proposition \ref{prop-hausdorff} reads as follows.

\begin{prop}\label{prop-hausdorff3}
The groupoid $\cG(\fX,\G,\Phi)$ is Hausdorff if and only if for all $[h]_x \in \cG(\fX, \G, \Phi)$ with $h \cdot x = x$, if there exists a sequence   $\{x_n\} \subset \fX$ which converges to $x$, and such that $[h]_{x_n} = [\id]_{x_n}$ for all $n$, then $[h]_{x} = [\id]_{x}$. 
\end{prop}

Taking the contrapositive of this statement, we obtain the following criterion:

\begin{criterion}\label{crit-nonH-element}
Let $(\fX,\G,\Phi)$ be a group action. Then $\cG(\fX,\G,\Phi)$ is a non-Hausdorff groupoid if and only if there exists a germ $[h]_x \in \cG(\fX, \G, \Phi)$ with $h \cdot x = x$, and a sequence   $\{x_n\} \subset \fX$ which converges to $x$ such that $[h]_{x_n} = [\id]_{x_n}$ for all $n$, and such that $[h]_{x} \ne [\id]_{x}$. 
 \end{criterion}
 
 We call a representative $h$ of such a germ a \emph{non-Hausdorff element} of $\cG(\fX,\G,\Phi)$. This is precisely the group element in Definition \ref{defn-nonH-rev}. 
 
 \subsection{Germinal groupoid of a locally quasi-analytic Cantor action}\label{sec-germ-LQA}
 
 From now on we assume that the action $(\fX,\G,\Phi)$ is minimal and equicontinuous.
 
 Suppose $(\fX,\G,\Phi)$ is locally quasi-analytic as in Definition \ref{def-LQA} with constant $\epsilon >0$. Then any $g \in \G$ which is trivial on an open subset of a set $U$ with $\diam(U)< \epsilon$ is trivial on $U$. Thus $(\fX,\G,\Phi)$ does not admit open sets $U$ of arbitrary small diameter such that $g|U$ is non-trivial, and $g|W = \id$ for some open $W \subset U$. Therefore, we have the following result, first observed in \cite[Proposition 2.5]{HL2018}.

\begin{criterion}\label{crit-hausdorff2}
If a minimal equicontinuous action $(\fX,\G,\Phi)$   is  locally quasi-analytic, then the germinal groupoid $\cG(\fX, \G, \Phi)$ is Hausdorff.     
 \end{criterion}
 
 The group $\G$ in Criterion \ref{crit-hausdorff2} can be chosen to be a countable or a profinite group.
 
 \begin{ex}
 {\rm
 Actions which are free or topologically free are locally quasi-analytic, see Section \ref{sec-LQA}, and so they have Hausdorff germinal groupoids. There are many examples of minimal equicontinuous actions of countable groups with this property. 
  For instance, in \cite{DHL2017} it was shown that any finite or separable profinite group can be realized as the discriminant group of a minimal equicontinuous action of a finite index torsion-free subgroup of ${\rm SL}(n,\mZ)$ with $n \geq 3$ sufficiently large. These actions are free, and so their germinal groupoids are Hausdorff.
  Essentially free actions on rooted trees in \cite[Section 5]{Grigorchuk2011} are topologically free, and so their associated germinal groupoids are Hausdorff. These actions include, among others, the action generated by the Bellaterra automaton, the action of the lamplighter group, the action of a solvable Baumslag-Solitar group $BS(1,3)$, and others.
  The minimal equicontinuous actions of surface groups in \cite{Joseph2023} are topologically free, and therefore their germinal groupoids are  Hausdorff.
 }
 \end{ex}

\subsection{Countable centralizer and the stabilizer direct limit groups}\label{sec-centrstab-def}

In this section we introduce new invariants of a minimal equicontinuous Cantor action $(\fX,\G,\Phi)$, the \emph{(direct limit) countable stabilizer and centralizer groups}. We use these invariants in Section \ref{sec-obstructiondiscrete} to develop an obstruction for $\G$ to have a non-Hausdorff element.

For a minimal equicontinuous action $(\fX,\G,\Phi)$, let $\cU$ be an adapted neighborhood system at $x \in \fX$ with group chain $\G^x_{\cU} = \{\G_\ell \mid \ell \geq 0\}$. Let $X_{\infty}$ denote the inverse limit space defined as in \eqref{eq-invlimit}, and $\whC_\infty$ and $\cD_x$ be the inverse limit representations of the profinite group $\fG(\Phi)$ and the stabilizer group $\fG(\Phi)_x$ respectively, see Theorem \ref{thm-quotientspace}. Recall that we assume that $\Phi$ is an effective action, so that $\G$ injects onto a dense subgroup of $\fG(\Phi)$, and so onto a dense subgroup of $\whC_\infty$, also denoted by $\G$.

As in Section \ref{sec-dirlimitgroups}, denote by $\widehat \G_\ell = \overline{\G}_\ell$ the topological closure of $\G_\ell$ in $\whC_\infty$, then $\{\widehat \G_\ell \mid \ell \geq 0\}$ is a decreasing chain of clopen subgroups of $\whC_\infty$. Let $K(\Phi) = \{K_\ell \mid \ell \geq 0\} \subset \cD_x$ and $Z(\Phi) = \{Z_\ell \mid \ell \geq 0\} \subset \cD_x$ be the systems of profinite stabilizer and centralizer subgroups respectively, as defined in \eqref{eq-Kchain} - \eqref{eq-Zchain}. Recall that $\Upsilon^x_s(\Phi)$ and $\Upsilon^x_c(\Phi)$ denote the direct limits of the systems $K(\Phi)$ and $Z(\Phi)$ with respect to the inclusion maps.

We are now interested in the subgroups of the countable group $\G$ with properties similar to those of subgroups in the chains $K(\Phi)$ and $Z(\Phi)$. 

The countable counterpart of the profinite discriminant group $\cD_x$ is the isotropy group of the action of $\G$ at $x$, which is the kernel of the group chain $\G^x_{\cU}$, namely,
$$\cK(\G^x_{\cU})  = \cD_x \cap \G = \bigcap_{\ell \geq 0} ~ \G_{\ell}.$$ 
Similarly to \eqref{eq-Kchain} - \eqref{eq-Zchain}, define the \emph{countable stabilizer} and \emph{centralizer} groups by
   \begin{align}\label{eq-KZG} K_\ell^\G & = \{ h \in \cK(\G^x_{\cU}) \mid [h,g] \in \cK(\G^x_{\cU})\textrm{ for all }g \in \G_\ell  \}, \\ \label{eq-KZGZ}
   Z_\ell^\G & = \{ h \in \cK(\G^x_{\cU}) \mid [h,g] = \id \textrm{ for all }g \in \G_\ell  \}.
   \end{align}

We now show that these groups can be obtained by simply intersecting $K_\ell$ and $Z_\ell$ with $\G$.

\begin{lemma}\label{lemma-intersection-simple}
For groups $K_\ell^G$ and $Z_\ell^G$ defined by \eqref{eq-KZG} - \eqref{eq-KZGZ}, and groups $K_\ell$ and $Z_\ell$ defined by \eqref{eq-Kchain} - \eqref{eq-Zchain}, we have $K_\ell^\G = K_\ell \cap \G $ and $Z_\ell^\G = Z_\ell \cap \G $. 
\end{lemma}

\proof We give a detailed proof for $K^\G_\ell$, the proof for $Z^\G_\ell$ is similar. We have
$$ K_\ell \cap \G  = \{ h \in \cK(\G^x_{\cU}) \mid \Ad(h)(\widehat g \, \cD_x) = \widehat g \, \cD_x \textrm{ for all }\widehat g \in \widehat \G_\ell  \}.$$

Note that the condition $\Ad(h)(\widehat g \, \cD_x) = \widehat g \, \cD_x$ for all $\widehat g \in \widehat \G_\ell$ is equivalent to the condition that $[h,\widehat g] \in \cD_x$ for all $\widehat g \in \G_\ell$. We first show that $K_\ell \cap \G \subset K^\G_\ell$. For that we must show that if $h \in \cK(\G^x_\cU)$ is such that $[h, \widehat g] \in \cD_x$ for all $\widehat g \in \widehat \G_\ell$, then $[h,g] \in \cK(\G^x_\cU)$ for all $g \in \G_\ell$. Since $\G_\ell \subset \widehat \G_\ell$ then we have 
  $[h,g] \in \cD_x$ for all $g \in \G_\ell$. Then since $h,g \in \G_\ell$, it follows that the commutator $[h,g] \in \G_\ell$, and so we must have $[h,g] \in \cD_x \cap \G_\ell = \cK(\G^x_\cU)$. Thus $h \in K^\G_\ell$.

Now lets us show that $K_\ell^\G \subset K_\ell \cap \G$. We have that $\whg  = (g_i) \in \widehat \G_\ell$, where $g_i \in \G_i \subset \G_\ell$ for $i \geq 0$.   Let $h \in K_\ell^\G$, then $ [h,\widehat g] = ([h,g_i])$. By assumption on $h$ we have $[h,g_i] \in \cK(\G^x_\cU)$ for each $i \geq 0$, and so $ [h, \widehat g] \cdot x = ([h,g_i]) \cdot x = x$. Then $[h, \widehat g] \in \cD_x$ and $h \in K_\ell \cap \G$. This shows that $K_\ell^\G = K_\ell \cap \G$. 
\endproof

\begin{remark}
{\rm
Although the countable group $\G$ is dense in the profinite group $\whC_\infty$, the kernel $\cK(\G^x_\cU)$ need not be dense in the discriminant group $\cD_x$, and the countable groups  $K^\G_\ell$ (resp. $Z^\G_\ell$) need not be dense in $K_\ell$ (resp. $Z_\ell$). For instance, \cite[Theorem 1.10]{DHL2017} proves that any finite or separable profinite group can be realized as the discriminant group of a stable minimal equicontinuous action of a finite index torsion-free subgroup of ${\rm SL}(n,\mZ)$ with $n \geq 3$ sufficiently large. These actions are free, so the kernel $\cK(\G^x_\cU)$ is trivial for any $x \in \fX$, while the discriminant group $\cD_x$ is a non-trivial finite or separable group. Since the action is stable, $K_\ell$ is trivial for $\ell$ sufficiently large. 

Another example is given by a family of minimal equicontinuous actions of the Heisenberg group $\cH$ constructed in \cite[Theorem 1.4]{HL2023}. Every action $(\fX,\cH,\Phi)$ in this family is topologically free, which implies that there exists $x \in \fX$ such that $\cK(\G^x_\cU) = \{\id\}$. Then for any $\ell \geq 0$ we have $K_\ell^\cH = \{\id\}$. On the other hand, $(\fX,\G,\Phi)$ is wild, which means that $K_\ell$ is non-trivial for all $\ell \geq 0$.

}
\end{remark}
 
We define the \emph{(direct limit) countable stabilizer and centralizer groups} by
 \begin{align}\label{eq-directlimk-countable} \Upsilon^{x,\G}_s(\Phi) & =\lim_{\longrightarrow} 
\cS(K_\ell^\G, \iota_\ell^{\ell+1},\mN) \subset \Upsilon^x_s(\Phi),  \\ \label{eq-directlimk-countable-Z} \Upsilon^{x,\G}_c(\Phi) & =\lim_{\longrightarrow} 
\cS(Z_\ell^\G, \iota_\ell^{\ell+1},\mN) \subset \Upsilon^x_c(\Phi).
 \end{align}
 
 We note the following property of countable stabilizer groups.
 
 \begin{lemma}
Let $(\fX,\G,\Phi)$ be a minimal equicontinuous Cantor actions. Then the action of $\G$ on $\fX$ is locally quasi-analytic if and only if the direct limit group $\Upsilon_s^\G(\Phi)$ of countable stabilizer subgroups is bounded.
\end{lemma}

\proof Let $\cU = \{U_\ell \mid \ell \geq 0\}$ be an adapted neighborhood basis at $x \in \fX$ with associated group chain $\G^x_\cU = \{\G_\ell \mid \ell \geq 0\}$, and the corresponding chain of normal cores $\cC = \{C_\ell \mid \ell \geq 0\}$, and let $K(\Phi) = \{K_\ell^\G \mid  \ell \geq 0\}$ and $\Upsilon_s^{x,\G}(\Phi) =\varinjlim \cS(K_{\ell}^\G, \iota^{\ell'}_{\ell} , \mN)$ be as above.
Suppose $\Upsilon_s^{x,\G}(\Phi)$ is unbounded, then there exists a cofinal set $\Lambda \subset \mN$ such that for $\ell' > \ell \in \Lambda$ the inclusion $\iota^{\ell'}_{\ell}: K^\G_\ell \to K^\G_{\ell'}$ is not an isomorphism. Then there exists $g \in \G$ such that $g|U_{\ell'} = \id$ and $g|U_\ell \ne \id$. Since the diameter of the sets $U_\ell$ tends to zero, this implies that $(\fX,\G,\Phi)$ is not locally quasi-analytic, see Definition \ref{def-LQA}.
\endproof

 While the profinite discriminant group $\cD_x$ is independent of the choice of $x \in \fX$ up to an isomorphism, a similar statement does not hold for the kernel $\cK(\G_\cU^x)$ of a group chain $\G^x_\cU$. It is easy to construct examples where $\cK(\G_\cU^x)$ is trivial for one choice of $x$, and non-trivial for another choice, see for instance \cite[Example 7.5]{DHL2016a}. However, since for $\ell \geq 0$ the elements in $K_\ell^\G$ act  trivially on clopen sets $U_\ell$, and not just at a point $x \in \fX$, their direct limit groups provide an invariant of conjugacy of group Cantor actions, as we show now. 
    
\begin{thm}\label{thm-isomupsilon-countable}
Let $(\fX,\G,\Phi)$ be a minimal equicontinuous action of a finitely generated group $\G$ on a Cantor set $\fX$. Then the direct limit isomorphism classes $\Upsilon_s^\G(\Phi)$ and $\Upsilon_c^\G(\Phi)$ of the groups $\Upsilon^{x,\G}_s(\Phi)$ and $\Upsilon^{x,\G}_c(\Phi)$ are invariants of the conjugacy class of the action $(\fX,\G,\Phi)$.
\end{thm}

\proof

Let  $(\fX,\G,\Phi)$ and  $(\fX',\G,\Psi)$ be conjugate minimal equicontinuous Cantor actions, that is, there is a homeomorphism $h \colon \fX' \to \fX$ such that for all $z \in \fX'$ and all $g \in \G$ we have $h(\Psi(g)(z)) = \Phi(g) (h(z))$. 
First we reduce the problem to considering two neighborhood bases for the \emph{same} action. 

Let  $\cU = \{U_{\ell} \subset \fX  \mid \ell \geq 0\}$ and $\cV = \{V_{\ell} \subset \fX' \mid \ell \geq 0\}$ be adapted neighborhood bases at $x \in \fX$ for the action $\Phi$, and at $z \in \fX'$ for the action $\Psi$ respectively. Then $\cU' = \{ U'_{\ell} = h(V_{\ell}) \subset \fX \mid \ell \geq 0\}$ is  an adapted neighborhood basis   at $y = h(z)$ for the action $\Phi$. If $\{H_\ell \mid \ell \geq 0\}$ is the group chain associated to $\cV$ and the action $\Psi$, and $\G^y_{\cU'} = \{\G_\ell' \mid \ell \geq 0\}$ is the group chain associated to $\cU'$ and the action $\Phi$, then $\G_\ell ' = H_\ell$ for all $\ell \geq 0$. 

To prove the theorem, we must show that the countable stabilizer and centralizer groups $\Upsilon_s^{x,\G}(\Phi)$ and $\Upsilon_c^{x,\G}(\Phi)$, associated  to the group chain $\G^x_{\cU}$  in $\G$ by \eqref{eq-KZG}, are isomorphic as direct limits to the countable stabilizer and centralizer groups $\Upsilon_s^{y,\G}(\Phi)$ and $\Upsilon_c^{y,\G}(\Phi)$ associated  to the group chain $\G^y_{\cU'}$.

  {\bf Part A (the same basepoint).} First assume that $x=y$, so we are given two  adapted neighborhood bases  at a common basepoint $x$,   $\cU = \{U_{\ell} \subset \fX  \mid \ell \geq 0\}$   and $\cU' = \{U'_{\ell} \subset \fX  \mid \ell \geq 0\}$, with group chains $\G^x_{\cU}$ and $\G^x_{\cU'}$.  
 
  As $\cU$ and $\cU'$ are both adapted neighborhood bases at $x$,  there exist   increasing sequences of indices   
 $1 \leq i_1 < i_2 < i _3 < \cdots$ and $1 \leq j_1 < j_2 < j_3 < \cdots$ such that we have a descending sequence of adapted clopen sets at $x$, where $i_0 = j_0 = 0$, 
 $$ \fX = U_0  = U'_0 \supset U_{i_1} \supset U'_{j_1} \supset U_{i_2} \supset U'_{j_2} \supset \cdots \ .$$
 Passing to a subsequence of coverings, we can   assume without loss of generality that $i_{\ell} = \ell$ and   $j_{\ell} = \ell$, for $\ell \geq 0$.
Introduce a common refinement $\cU'' = \{U''_{\ell} \mid \ell \geq 0\}$ of these chains of clopen sets, where $U''_{2\ell} = U_{\ell}$ and $U''_{2\ell -1} = U'_{\ell}$.

 Let   $\G^x_{\cU''} = \{\G''_{\ell}\mid \ell \geq 0\}$  be the group chain  associated to $\cU''$, then $\G''_{2\ell} = \G_{\ell}$ and $\G''_{2\ell -1} = \G'_{\ell}$. Let $X_{\infty}$, $X'_{\infty}$ and $X''_{\infty}$ denote the inverse limit spaces  defined as in \eqref{eq-invlimit} by the group chains   $\G^x_{\cU}$,  $\G^x_{\cU'}$ and  $\G^x_{\cU''}$, respectively.

 By  the discussion in Section~\ref{sec-adapted},  there are homeomorphisms 
 $f_x \colon \fX \to X_{\infty}$, $f_x' \colon \fX \to X'_{\infty}$ and $f_x'' \colon \fX \to X''_{\infty}$ which intertwine the $\G$-actions on these spaces.
Introduce the basepoint preserving homeomorphisms  
  \begin{align*}\tau & =  f_x \circ (f_x'')^{-1} \colon X''_{\infty} \to X_{\infty}: (g''_{\ell} \, \G_\ell'') \mapsto (g''_{2\ell} \G_\ell) , \\
   \tau' & =  f_x' \circ (f_x'')^{-1} \colon X''_{\infty} \to X'_{\infty}: (g''_{\ell} \G_\ell'') \mapsto (g''_{2\ell-1}\G_\ell').\end{align*}
Here we use the full notation for elements of $X_\infty$ (resp. $X_\infty'$, $X_\infty''$) as sequences of cosets $(g_\ell \G_\ell)$ (resp. $(g_\ell' \G_\ell')$, $(g_\ell'' \G_\ell'')$), instead of the short notation $(g_\ell)$ (resp. $(g_\ell')$, $(g_\ell'')$) which we use in the rest of the paper, to make the definitions of $\tau$ and $\tau'$ clearer.
 
 Recall from Section \ref{sec-profinite} that, given $\G^x_\cU = \{\G_\ell \mid \ell \geq 0\}$, there is a descending chain $\cC_\cU = \{C_\ell \mid \ell \geq 0\}$, where $C_\ell$ is the normal core of $\G_\ell$. Similarly, there are chains $\cC_{\cU'} = \{C_\ell' \mid \ell \geq 0\}$ and $\cC_{\cU''} = \{C_\ell'' \mid \ell \geq 0\}$ for the chains of normal cores of the groups in $\G^x_{\cU'}$ and $\G^x_{\cU''}$ respectively.
    
 Let $\whC_{\infty}$, $\whC'_{\infty}$ and $\whC''_{\infty}$ denote the inverse limit groups  defined as in \eqref{eq-Cspace}  by the chains  $\cC_{\cU}$,  $\cC_{\cU'}$ and  $\cC_{\cU''}$, respectively. 
By  Theorem~\ref{thm-quotientspace}, there are topological isomorphisms 
\begin{align}\label{eq-cisom}\widehat f \colon \fG(\Phi) \to \whC_{\infty} ~ , ~ \widehat f' \colon \fG(\Phi) \to \whC'_{\infty} ~ , ~ \widehat f'' \colon \fG(\Phi) \to \whC''_{\infty} \ ,\end{align}
which map $\G \subset \fG(\Phi)$ onto dense subgroups of $\whC_\infty$, $\whC_\infty'$ and $\whC_\infty''$. The compositions of the isomorphisms in \eqref{eq-cisom} give topological group isomorphisms  
\begin{align*}\whtau & = \widehat f \circ (\widehat f'')^{-1} \colon \whC''_{\infty} \to \whC_{\infty} : (g''_{\ell}C_\ell'') \mapsto (g''_{2\ell}C_\ell),\\
   \whtau' & = \widehat f' \circ (\widehat f '')^{-1} \colon \whC''_{\infty} \to \whC'_{\infty}: (g''_{\ell} C_\ell'') \mapsto (g''_{2\ell-1}C_\ell') .\end{align*}
The maps $\whtau$ and $\whtau'$ map the dense subgroup $\widehat f''(\G)$ of $\whC_\infty''$ onto the dense subgroups $\widehat f(\G)$ and $\widehat f'(\G)$ of $\whC_\infty$ and $\whC_\infty'$ respectively. This is a key property of the construction, which need not hold for the case of distinct basepoints in {\bf Part B}.

Denote by $\cD_x$, $\cD_x'$ and $\cD_x''$ the images of the isotropy group $\fG(\Phi)_x$ under the maps in \eqref{eq-cisom} respectively. 
Then 
   the restrictions    $\whtau \colon \cD''_x \to  \cD_x  \subset \whC_{\infty}$ and $\whtau' \colon \cD''_x \to  \cD'_x \subset \whC'_{\infty}$ are isomorphisms of the discriminant groups. Then for the kernels of the group chains we have
     $$\cK(\G^x_\cU) = \cD_x \cap \widehat f(\G), \quad \cK(\G^x_{\cU'}) = \cD_x' \cap \widehat f'(\G), \quad \cK(\G^x_{\cU''}) = \cD_x'' \cap \widehat f''(\G).$$
  Since $\whtau$ and $\whtau'$ preserve the discriminant groups, and map the dense subgroup $\widehat f''(\G)$ onto the dense subgroups $\widehat f(\G)$ and $\widehat f'(\G)$ respectively, they preserve the kernels of the group chains. Namely, they restrict to the isomorphisms
    $$\whtau: \cK(\G^x_{\cU''}) \to \cK(\G^x_{\cU}), \quad \textrm{and} \quad \whtau': \cK(\G^x_{\cU''}) \to \cK(\G^x_{\cU'}).$$
Now consider the countable stabilizer groups, defined in Lemma \ref{lemma-intersection-simple},
\begin{eqnarray*}
K_{\ell}^\G   & = &   \{h \in  \cK(\G^x_\cU) \mid [h,\widehat g] \in \cD_x \ {\rm for ~all}~ \whg \in \widehat \G_{\ell} \} \subset \cD_x \subset \widehat \G_\ell,  \\
(K'')_{\ell}^\G  & = &  \{h \in  \cK(\G^x_{\cU''}) \mid [h,\widehat g] \in \cD_x'' \ {\rm for ~all}~ \whg \in \widehat \G_{\ell}'' \}   \subset \cD_x'' \subset \widehat \G_\ell''.
\end{eqnarray*}
Here $\widehat \G_\ell$ is the closure of the subgroup $\widehat f (\G_\ell)$ in $\whC_\infty$, where $\G_\ell \in \G^x_\cU$, and $\widehat \G_\ell''$ is the closure of the subgroup $\widehat f'' (\G_\ell'')$ in $\whC_\infty''$, where $\G_\ell'' \in \G^x_{\cU''}$. Note that for $\ell \geq 0$ there is a group isomorphism $\whtau: \widehat f'' (\G_{2\ell}'') \to \widehat f (\G_\ell)$, and so there is a group  isomorphism $\whtau \colon \widehat \G''_{2\ell} \to \widehat \G_{\ell}$. Since $\whtau$ preserves the discriminant group and the kernel of the group chain, it restricts to a group isomorphism $\whtau: (K'')_{2 \ell}^\G \to K_{\ell}^\G$, for any $\ell \geq 0$. Since the set $\Lambda = \{2 \ell \mid  \ell \geq 0\}$ is cofinal in $\mN$, by Propositions~\ref{prop-directedgroupisomorphism} and \ref{prop-directedgroupcofinal} there is an isomorphism of direct limit groups
\begin{equation}\label{eq-isoisotropy}
 \underrightarrow{\whtau} \colon (\Upsilon'')_s^{x,\G} =  \varinjlim \cS((K'')_{\ell}^\G, \iota^{\ell'}_{\ell} , \mN) \to \Upsilon_s^{x,\G} = \varinjlim \cS(K_{\ell}^\G, \iota^{\ell'}_{\ell} , \mN) \ .
\end{equation}
Similarly we can show that the direct limit groups $(\Upsilon'')_s^{x,\G} =\varinjlim \cS(K_{\ell}, \phi^{\ell'}_{\ell} , \mN)$ and $(\Upsilon')_s^{x,\G} = \varinjlim \cS(K'_{\ell}, \psi^{\ell'}_{\ell} , \mN)$ of countable stabilizer subgroups are  isomorphic. Therefore, $\Upsilon_s^{x,\G}$ and $(\Upsilon')_s^{x,\G}$ are isomorphic, representing the same isomorphism class $\Upsilon_s^\G$.
  
 The proof that the direct limit groups $\Upsilon_c^{x,\G}$ and $(\Upsilon')_c^{x,\G}$ are isomorphic is similar, using the same adapted neighborhood system $\cU''$ as above, so we omit it.     

{\bf Part B (distinct basepoints).} Next,  consider the case where $x \ne y$, and we are given    adapted neighborhood bases    $\cU = \{U_{\ell} \subset \fX  \mid \ell \geq 0\}$  at $x$  and $\cU' = \{U'_{\ell} \subset \fX  \mid \ell \geq 0\}$ at $y$, with corresponding group chains $\G^x_{\cU} = \{\G_{\ell} \mid \ell \geq 0\}$ and $\G^y_{\cU'} = \{\G'_{\ell}\mid \ell \geq 0\}$.  

Since $\diam(U_\ell')$ tend to zero with $\ell$, and possibly restricting to a subsequence, for any $\ell \geq 0$ we can choose $g_\ell \in \G$ such that $g_\ell \cdot U_\ell' \owns x$, and $g_{\ell+1} \cdot U_{\ell+1}' \subset g_\ell \cdot U_\ell'$. Then $\cU'' = \{U_\ell'' = g_\ell \cdot U'_{\ell} \subset \fX  \mid \ell \geq 0\}$  is an  adapted neighborhood basis  at $x$, with the associated group chain   
 $\G^x_{\cU''} = \{\G''_{\ell} = g_{\ell} \ \G'_{\ell} \ g_{\ell}^{-1} \}$. 
 
 Let $X'_{\infty}$ and $X''_{\infty}$ denote the inverse limit spaces  defined as in \eqref{eq-invlimit} by the group chains   $\G^y_{\cU'}$ and  $\G^x_{\cU''}$, respectively. By  the discussion in Section~\ref{sec-adapted},  there is a homeomorphism
 $f_y' \colon \fX \to X'_{\infty}$, associated to the group chain $\G^y_{\cU'}$, and a homeomorphism $f_x'' \colon \fX \to X''_{\infty}$, associated to the group chain $\G^x_{\cU''}$, which intertwine the $\G$-actions on the corresponding spaces. 
  
 Let $C_\ell'$ be the normal core of $\G_\ell'$, then $C_\ell'$ is also the normal core of $\G_\ell''$. Let $\cC_{\cU'} = \{C_\ell' \mid \ell \geq 0\}$ be the group chain of normal cores, and denote by $\whC'_{\infty}$ the inverse limit group  defined as in \eqref{eq-Cspace}  by the chain  $\cC_{\cU'}$.
By  Theorem~\ref{thm-quotientspace}, there is topological isomorphism 
\begin{align}\label{eq-cisom-conj} \widehat f' \colon \fG(\Phi) \to \whC'_{\infty} ~ ,\end{align}
which maps $\G$ onto a dense subgroup $\widehat f'(\G)$ of $\whC_{\infty}'$. Denote by $\cD_y = f'(\fG(\Phi)_y)$ and $\cD_x'' = f'(\fG(\Phi)_x)$ the images of the isotropy groups of the action of $\fG(\Phi)$ at $y$ and $x$ respectively. We note that $\fG(\Phi)_x$ and $\fG(\Phi)_y$ are conjugate but not necessarily equal, and a similar statement holds for $\cD_x''$ and $\cD_y'$.

 Conjugation by $g_\ell$ induces bijections of finite cosets spaces
  $$\tau_\ell: X_\ell' \to X_\ell'': g \, \G_\ell' \to (g_\ell g g_\ell^{-1})\G_\ell'',$$ 
  and so induces a homemorphism of the inverse limits of coset spaces
  \begin{align*}\tau' &: X_\infty' \to X_\infty'': (h_\ell \G_\ell' ) \mapsto (g_\ell h_\ell g_\ell^{-1} \G_\ell''). \end{align*}
Namely, $\tau'$ is induced by conjugation by $\whg = (g_\ell) : = (g_\ell C_\ell) \in \whC_\infty'$, and the conjugation also induces the map of profinite groups
  $\whtau': \whC_\infty' \to \whC_\infty'$. Then 
   $$\whtau'(\cD_y')= \whg \, \cD_y' \whg^{-1} = \cD_x''.$$ 
However, if $h \in \cD_y' \cap \widehat f'(\G)$, the conjugate $\whg  \, h \, \whg^{-1} \in \cD_x''$ need not be in $\widehat f'(\G)$. This is the subtlety which arises in working with countable subgroups in profinite groups. In particular, for the kernels of the group chains $\cK(\G^y_{\cU'}) = \cD_y' \cap \widehat f'(\G)$ and $\cK(\G^x_{\cU''}) = \cD_x'' \cap \widehat f'(\G)$, the image $\whtau'(\cK(\G^y_{\cU'}))$ need not have $\cK(\G^x_{\cU''})$ as its image.

We now show that the direct limit countable stabilizer and centralizer groups associated to $\G^y_{\cU'}$ and $\G^x_{\cU''}$ are nevertheless isomorphic. 

Define the countable stabilizer groups, as in \eqref{eq-KZG}, by
\begin{eqnarray*}
(K')_{\ell}^\G   & = &   \{h \in  \cK(\G^y_{\cU'}) \mid [h,g] \in \cK(\G^y_{\cU'}) \ {\rm for ~all}~ g \in \G_{\ell}' \} \subset \cD_y' \subset \widehat \G_\ell',  \\
(K'')_{\ell}^\G    & = &  \{h \in  \cK(\G^x_{\cU''}) \mid [h, g] \in \cK(\G^y_{\cU''})  \ {\rm for ~all}~ g \in \G_{\ell}'' \}   \subset \cD_x'' \subset \widehat \G_\ell''.
\end{eqnarray*}

Denote $g^h = h gh^{-1}$. Note that the restriction $g|{U_\ell'} = \id$ if and only if $g^{g_\ell}|U_\ell'' = \id$.

Let $h \in (K')_{\ell}^\G$, and recall that $[h,g] \in \cK(\G^y_{\cU'})$ for all $g \in \G_\ell'$ implies that $h|U_\ell' = \id$, see Section \ref{sec-dirlimitgroups}. Since $h,g_\ell \in \G_\ell$, we have $h^{g_\ell} \in \G_\ell$, and by above remarks $h^{g_\ell}|U_\ell'' = \id$. Since $x \in U_\ell''$, then $h^{g_\ell} \in (K'')_\ell^\G$ even if $x \ne g_\ell \cdot y$. Thus the restriction $\whtau': (K')^\G_\ell \to \cD_x'' \subset \widehat \G_\ell''$ is an isomorphism onto a subgroup of $(K'')_\ell^\G$. Applying a similar argument to the elements in $(K'')_\ell^\G$ and conjugating by $g_\ell^{-1}$ we obtain that the restriction $(\whtau')^{-1}: (K'')^\G_\ell \to \cD_y' \subset \widehat \G_\ell'$ is an isomorphism onto a subgroup of $(K')_\ell^\G$. 
The composition $(\whtau')^{-1} \circ \whtau': (K')^\G_\ell \to (K')^\G_\ell$ is the identity map, which implies that $\widehat \tau': (K')_\ell^\G \to (K'')_\ell^\G$ is surjective and so an isomorphism onto $(K'')_\ell^\G$. Then by Proposition~\ref{prop-directedgroupisomorphism} there is an isomorphism of direct limit groups
\begin{equation}\label{eq-isoisotropy-1}
 \underrightarrow{\whtau} \colon  (\Upsilon')_s^{x,\G} =  \varinjlim \cS((K')_{\ell}^\G, \iota^{\ell'}_{\ell} , \mN) \to  (\Upsilon'')_s^{x,\G} = \varinjlim \cS((K'')_{\ell}^\G, \iota^{\ell'}_{\ell} , \mN) \ .
\end{equation}
Next, $\cU$ and $\cU''$ are adapted neighborhood centered at the same point $x \in \fX$, so we use the argument in {\bf Part A} to obtain an isomorphism $(\Upsilon'')_s^{x,\G} \to \Upsilon_s^{x,\G}$. Composing it with \eqref{eq-isoisotropy-1} we obtain an isomorphism $(\Upsilon')_s^{y,\G} \to \Upsilon_s^{x,\G}$ of the direct limit groups of countable stabilizer subgroups associated to the neighborhood systems $\cU$ and $\cU''$. Thus the isomorphism class $\Upsilon^\G_s(\Phi)$ is invariant under conjugacy of minimal equicontinuous actions.

The proof for the direct limit countable centralizer group is similar, so we omit it.

\endproof

In Definition \ref{defn-types-wild-countable} we introduce a classification of actions of countable groups, similar to the classification which involves the actions of their completions.

\begin{defn}\label{defn-types-wild-countable}
Let   $(\fX,\G,\Phi)$ be a minimal equicontinuous Cantor action, and suppose the action of $\G$ on $\fX$ is not locally quasi-analytic. Then $(\fX,\G,\Phi)$ is:
\begin{enumerate}
\item \emph{countably wild of finite type}  if    the countable stabilizer  group $\Upsilon_s^\G(\Phi)$   is unbounded, and represented by a chain of finite groups;
\item \emph{countably wild of flat type}  if    the countable stabilizer  group $\Upsilon_s^\G(\Phi)$   is unbounded, and $\Upsilon_c^\G(\Phi) = \Upsilon_s^\G(\Phi)$; 
\item \emph{countably dynamically wild}  if    the countable stabilizer  group $\Upsilon_s^\G(\Phi)$   is unbounded, and is not of flat type.
\end{enumerate}
\end{defn} 

The relationship between the notions in Definition \ref{defn-types-wild-countable} and \ref{defn-types-wild} is given by the following corollary of Theorem \ref{thm-isomupsilon-countable}.

\begin{cor}\label{cor-incl-prof-count}
Let $(\fX,\G,\Phi)$ be a minimal equicontinuous action, and suppose that the action of $\G$ is wild. Then: 
\begin{enumerate}
\item If  $(\fX,\G,\Phi)$ is wild of finite type, then $ (\fX,\G,\Phi)$ is either locally quasi-analytic, or countably wild of finite type.
\item If  $(\fX,\G,\Phi)$ is wild of flat type, then $ (\fX,\G,\Phi)$ is either locally quasi-analytic, or countably wild of flat type.
\item If  $(\fX,\G,\Phi)$ is countably dynamically wild, then $ (\fX,\G,\Phi)$ is dynamically wild.
\end{enumerate}
\end{cor}

\proof
For the proof of all statements it is useful to note the following simple fact: for all $\ell \geq 1$, the mappings $Z_\ell \to K_\ell$ are inclusions of subgroups, and thus the intersection with the countable subgroup $\G$ is preserved under this mapping. Namely, if $g \in Z^\G_\ell = Z_\ell \cap \G$, then $g \in K^\G_\ell = K_\ell \cap \G$. Thus if the inclusion $Z_\ell^\G \to K_\ell^\G$ is not surjective, then the inclusion $Z_\ell \to K_\ell$ is not surjective which implies (3).

For statements (1) and (2) we note that it is possible to have an action $\G$ which is locally quasi-analytic (LQA), or even topologically free, with wild closure, see \cite[Theorem 1.4]{HL2023}. So if $(\fX,\G,\Phi)$ is wild, then the action of $\G$ on $\fX$ may be either LQA or not LQA. If $(\fX,\G,\Phi)$ is not LQA and it is wild of finite type, then the stabizer group $\Upsilon_s(\Phi)$ is the direct limit of finite groups, and therefore the group $\Upsilon_s^\G(\Phi)$ is also the direct limit of finite groups. Then $(\fX,\G,\Phi)$ is countably wild of finite type.

If $(\fX,\G,\Phi)$ is not LQA and it is wild of flat type, then the inclusion $\Upsilon_c(\Phi) \to \Upsilon_s(\Phi)$ is an isomorphism. Since the elements of $\G$ are mapped onto elements of $\G$ under inclusions, this implies that $\Upsilon_c^\G(\Phi) \to \Upsilon_s^\G(\Phi)$ is an isomorphism. Then $(\fX,\G,\Phi)$ is countably wild of flat type.
\endproof

\subsection{An obstruction to the existence of a non-Hausdorff element}\label{sec-obstructiondiscrete}

We now show that if the germinal groupoid $\cG(\fX,\G,\Phi)$ is non-Hausdorff, then the action $(\fX,\G,\Phi)$ must be countably dynamically wild.

Recall that for the direct limits $\Upsilon_s(\Phi) $ and $\Upsilon_c(\Phi) $ of profinite groups we have the following criterion.

\begin{thm} \cite[Theorem 1.7]{HL2021}\label{thm-incl-profinite}
Let $(\fX,\G,\Phi)$ be a minimal equicontinuous Cantor action. If the action of the closure $\fG(\Phi) = \overline{\Phi(\G)}$ on $\fX$ has a non-Hausdorff element, then $\Upsilon_c(\Phi) \to \Upsilon_s(\Phi)$ is a proper inclusion.
\end{thm}

Theorem \ref{thm-incl-profinite} can be used as an obstruction to the existence of a non-Hausdorff element in the group $\G$: if the inclusion $\Upsilon_c^x(\Phi) \to \Upsilon_s^x(\Phi)$ is a group isomorphism, then by Theorem \ref{thm-incl-profinite} the profinite group $\fG(\Phi)$, and therefore its dense countable subgroup $\G \subset \fG(\Phi)$ does not contain non-Hausdorff elements. 

However, it is conceivable that $\fG(\Phi)$ contains non-Hausdorff elements, while $\G$ does not, and Theorem \ref{thm-incl-profinite} does not distinguish between these two cases. In the rest of this section, we develop a more sensitive criterion which allows to rule out the existence of non-Hausdorff elements in $\G$, using the countable stabilizer and centralizer groups defined in Section \ref{sec-centrstab-def}.

We now prove our main technical result. The key point in the proof is that we have to work with a countable group acting on $\fX$ minimally but not transitively, and so   special care must be taken when making certain choices.

\begin{thm} \label{thm-incl-countable}
Let $(\fX,\G,\Phi)$ be a minimal equicontinuous Cantor action. If $\G$ has a non-Hausdorff element, then $(\fX,\G,\Phi)$ is countably dynamically wild, i.e. the inclusion $\Upsilon_c^{\G}(\Phi) \to \Upsilon_s^{\G}(\Phi)$ is proper.
\end{thm}

\proof  

Let $\cU = \{U_{\ell} \subset \fX  \mid \ell \geq 0\}$ be an adapted neighborhood basis for the action $\Phi$ at $x$, and let  $\G^x_{\cU} = \{\G_{\ell}\mid \ell \geq 0 \}$ be the associated group chain.  Denote by $X_{\infty}$ the inverse limit space  defined as in \eqref{eq-invlimit}, and recall that there is a homeomorphism $f: \fX \to X_\infty$ which intertwines the actions of $\G$ on the corresponding spaces. Denote  $\widetilde U_\ell = f(U_\ell)$ for each $\ell \geq 0$.

Let $\cC_{\cU}$ be the group chain of normal cores of the groups in $\G^x_\cU$, and denote by $\whC_{\infty}$ the inverse limit group  defined as in \eqref{eq-Cspace}. By  Theorem~\ref{thm-quotientspace}, there is a topological isomorphism 
$\widehat f \colon \fG(\Phi) \to \whC_{\infty}$, and we denote also by $\G$ the image $\widehat f(\G)$, and set $\cD_x = \widehat f(\fG(\Phi)_x)$. Denote by $\cK(\G^x_\cU) = \bigcap_{ \ell \geq 0} \G_\ell$ the kernel of $\G^x_\cU$. For $\ell \geq 0$, the topological closure $\widehat \G_\ell = \overline{\G_\ell}$ is a clopen subgroup  of $\whC_\infty$, and we have a descending chain $\{\widehat \G_\ell \mid \ell \geq 0\}$ of clopen subgroups such that $\cD_x = \bigcap_{\ell \geq 0} \widehat \G_\ell$. 

For each $\ell \geq 0$, denote by $\widehat f_\ell: \whC_\infty \to \G/C_\ell$ the projection. 
By \cite[Proposition 5.7]{DHL2016a} the adapted neighborhood basis $\cU$ can be chosen so that the associated group chain $\G^x_\cU$ is in the \emph{normal form}, which means that the restriction $\widehat f_\ell(\cD_x) \to \G_\ell/C_\ell$, and so each $\G_i/C_i \to \G_\ell/C_\ell$ for $\ell \geq 0$, is surjective. We assume that $\cU$ and $\G^x_\cU$ have this property.

For the normal core $C_\ell \subset \G_\ell$, denote by $\whC_\ell = \overline{C_\ell} \subset \whC_\infty$ the topological closure of $C_\ell$. Since $C_\ell$ has finite index in $\G$, then $\whC_\ell$ is a clopen normal subgroup of $\whC_\infty$.  We will need the following result for our proof.

\begin{lemma}\label{lemma-normalform-transitive}
If the group chain $\G^x_\cU$ is in the normal form, then $\whC_\ell \, \cD_x  = \widehat \G_\ell$, and the following holds: 
\begin{enumerate}
\item The profinite group $\whC_\ell$ acts transitively on the clopen set $\widetilde U_\ell$, and
\item The countable group $C_\ell$ acts minimally on the clopen set $\widetilde U_\ell$.
\end{enumerate}
\end{lemma}

\proof
The inclusion $\whC_\ell \, \cD_x  \subseteq \widehat \G_\ell$ holds since  $\whC_\ell, \cD_x \subset \widehat \G_\ell$.  We have to show the reverse inclusion, $\widehat \G_\ell \subseteq \whC_\ell \, \cD_x$. By Theorem \ref{thm-quotientspace} we can write 
\begin{align*}\whC_\ell &= \{ (h_i C_i )\mid i \geq 0, h_i \in C_\ell\}, && \cD_x  = \{(g_i C_i) \mid i \geq 0, g_i \in \G_i\}. 
\end{align*}

For each $i \geq 0$, the subgroup $C_i$ is normal, and so its left and right cosets in $\G_\ell$ are equal. Then
  $$\whC_\ell \, \cD_x = \{(h_i g_i C_i) \mid i \geq 0, h_i \in C_\ell,\, g_i \in \G_i\}.$$
We will show that the countable subgroup $\G_\ell$ of $\G$ is dense in $\whC_\ell\, \cD_x$, and therefore we have $\widehat \G_\ell  = \widehat C_\ell \, \cD_x$ for its closure. For that we show that the restrictions $\widehat f_i \circ \widehat f:\G_\ell \to \G/C_i : g \mapsto g C_i$, where $\G_\ell \subset \fG(\Phi)$, are surjective onto the image $\widehat f_i (\whC_\ell \, \cD_x)$, for each $i \geq \ell$. Then by \cite[Lemma 1.1.7]{RZ2000} the group $\G_\ell$ maps onto a dense subgroup of $\widehat C_\ell \, \cD_x$. Thus for each $g \in \G_\ell$ and each $i \geq \ell$ we must find $h_i \in C_\ell$ and $g_i \in \G_i$ such that $g C_i = h_ig_i C_i$. Note that this clearly holds for $i = \ell$. 

Consider the action of $C_\ell$, $\G_\ell$ and $\G_i$ on $\G_\ell/C_i$, for $i \geq \ell$. Let $k = |G_\ell: C_\ell|$, the index of $C_\ell$ in $\G_\ell$. Let $g_1,\ldots, g_k$ be the representatives of the cosets of $C_\ell$ in $\G_\ell$. 

 Since $C_i \subset C_\ell$, then there is the inclusion of cosets $\G_\ell/C_i \to \G_\ell/C_\ell$, and we can partition the finite set $\G_i/C_i$ into the sets $\mathcal P = \{P_1,\ldots,P_k\}$, so that  $hC_i \in P_j$ if and only if $hC_i \subset g_j C_\ell$. Since $C_\ell$ is a normal subgroup, we have for each $h \in C_\ell$ and each $g_j \in \G_\ell$,
  $$h \cdot g_j C_\ell = g_j h'C_\ell = g_j C_\ell,$$
  for some $h'\in C_\ell$. Thus the sets of the partition $\mathcal P$ are preserved by the action of $C_\ell$. 
  
  Moreover, the action of $C_\ell$ is transitive on the cosets of $C_i$ in each $P_j$. Indeed, let $g_1C_i, g_2 C_i \in P_j$, where $g_1,g_2 \in \G_\ell$. By definition of $P_j$ these cosets are in the same left coset of $\G_\ell/C_\ell$. Since $C_\ell$ is a normal subgroup of $\G_\ell$, its left and right cosets are equal, and $g_1,g_2$ are in the same right coset of $C_\ell$. Then $g_1 g_2^{-1} \in C_\ell$. The left action of $g_1g_2^{-1}$ takes $g_2 C_\ell$ to $g_1 C_\ell$, which shows that the action of $C_\ell$ on the cosets in $P_j$ is transitive.
  
Now let $g \in \G_\ell$. Since the group chain $\G^x_\cU$ is in the normal form, the inclusions $\G_i/ C_i \to \G_\ell /C_\ell$ are surjective for all $i \geq \ell$, and so there exists $g_i \in \G_i$ such that $g_i C_i \in g C_\ell$. In particular, this means that $g_iC_i$ and $gC_i$ are in the same set $P_j$ of the partition $\cP$, and so there exists $h_i \in C_\ell$ such that $h_i g_iC_i = gC_i$. Thus the map $\widehat f_i \circ \widehat f:\G_\ell \to \widehat f_i(\whC_\ell \, \cD_x)$ is surjective. This argument holds for all $i \geq \ell$, thus the image of the map $\widehat f: \G_\ell \to \whC_\infty \, \cD_x$ is dense in $\widehat C_\ell \cD_x$, which implies that $\widehat \G_\ell = \widehat C_\ell \, \cD_x$.

We have shown that $\widehat \G_\ell = \widehat C_\ell \, \cD_x$. Then using the homogeneous model for the action, we have that
\begin{align}\label{whu-whc}\widetilde U_{\ell} = \widehat \G_{\ell}/ \cD_x = (\whC_{\ell} \ \cD_x)/ \cD_x = \whC_{\ell}/(\whC_{\ell} \cap \cD_x) \ .\end{align}
To see that the last equality is true, consider the map $\varphi:\whC_\ell \to \whC_\ell \, \cD_x / \cD_x: \whh \mapsto \whh \, \cD_x$. This map is clearly surjective, since, for any $\whh \in \whC_\ell$ and $\whd \in  \cD_x$ we have $\whh \, \whd \, \cD_x  = \whh \, \cD_x =  \varphi(\whh)$. If $\varphi(\whh_1) = \varphi(\whh_2)$, then $\whh_1^{-1}  \whh_2 \in \cD_x$, and also $\whh_1^{-1} \whh_2 \in \whC_\ell$, since $\whC_\ell$ is a group. Thus there is a bijection of coset spaces $\whC_\ell /(\whC_\ell \cap \cD_x) \to \whC_\ell \, \cD_x / \cD_x$.

In particular, \eqref{whu-whc} implies that $\whC_{\ell}$ acts transitively on $\widetilde U_\ell$. Since $C_\ell$ is dense in $\whC_\ell$, then $C_\ell$ acts minimally on $\widetilde U_\ell$.
   \endproof  

Next, define the chain
$\{K_\ell^\G  \mid \ell \geq 0\}$ of countable stabilizer groups, and    the chain   $\{Z_\ell^\G  \mid \ell \geq 0\}$   of countable centralizer groups as in \eqref{eq-KZG}-\eqref{eq-KZGZ}.

Let $g \in \G$ be a non-Hausdorff element at $x \in \fX$. Then the germ of $g$ at $x$ is non-trivial, and there exists a sequence   $\{x_i \mid i \geq 1\} \subset \fX$   of distinct points converging to $x$, with  $g \cdot x_i = x_i$ for all $i \geq 1$, and  clopen subsets $x_i \in W_i \subset  \fX$ such that  the restriction of $g$   to $W_i$ is the identity. Denote also by $x$ and $x_i$ the images of $x$ and $x_i$, $i \geq 1$, under $f$ in $X_\infty$, and set $\widetilde W_i = f(W_i)$.  

Note that $g \notin K_\ell^\G$ for any $\ell \geq 0$, since by assumption $g$ is non-trivial on any $\widetilde U_\ell$. There exists an $\ell_0$ such that for each $\ell \geq \ell_0$ the group $K_\ell^\G$ is non-trivial. Indeed, by minimality for each $i \geq 0$ there exists $g_i \in \G$ such that $g_i \cdot x \in \widetilde W_i$. Choose $k_i \geq i$ large enough so that $g_i(\widetilde U_{k_i}) \subset \widetilde W_i$. Then 
  $h_i = g_i^{-1} g g_i \in K_{k_i}^\G$, and for $\ell > k_i$ the group $K_\ell^\G$ is non-trivial.
  
 By increasing $k_i$ in the previous paragraph, if necessary, we can arrange that $\widetilde U_{k_i} \cap \widetilde W_i = \emptyset$, and then  $g_i \cdot \widehat U_{k_i} \cap \widehat U_{k_i} = \emptyset$. Moreover, as $x_i$ limits to $x$ there exists $j > k_i$ such that $x_j \in \widetilde U_{k_i}$ and so $x_j \not\in \widehat W_i$. Thus, passing to   subsequences chosen recursively, 
we can assume that we have a sequence $\{x_\ell\} \to x$, an adapted neighborhood system $\widetilde \cU = \{\widetilde U_\ell \mid \ell \geq 0\}$, and a collection of clopen sets $\{\widetilde W_\ell \mid \ell \geq 0\}$ with the following properties: 
\begin{equation}\label{eq-indexing}
 \widetilde U_{\ell} \cap \widetilde W_{\ell} = \emptyset    ~ , ~  x_{\ell} = g_{\ell} \cdot x \in \widetilde W_\ell ~\textrm{for} ~g_\ell \in \G,   ~    ~ g_{\ell} \cdot \widetilde U_{\ell} \subset \widetilde W_{\ell}    ~ , ~  x_{\ell +1} \not\in \widetilde W_{\ell} ~ .
\end{equation}
Then for each $\ell \geq 0$ we have that $h_\ell = g_\ell^{-1} g g_\ell \in K_\ell^\G$.
 
 Set  $y_{\ell} = g_{\ell}^{-1} \cdot x$ and observe that $h_{\ell} \cdot y_{\ell} = y_{\ell}$ and that $y_{\ell} \not\in \widetilde U_{\ell+1}$ as $g_{\ell} \cdot \widetilde U_{\ell+1} \cap \widetilde U_{\ell+1} = \emptyset$. Also, $h_\ell$ is not the identity on any open neighborhood of $y_\ell$, and so $h_\ell$ has a non-trivial germ at $y_\ell$. Moreover, $h_\ell$ is non-Hausdorff at $y_\ell$, and so there is $W' \subset g_\ell^{-1} \cdot \widetilde U_\ell$ such that $h_\ell| W'= \id$.
 
 Since $C_\ell$ is a normal subgroup of $\G$, its action preserves $\widetilde U_\ell$ and $g_\ell^{-1}(\widetilde U_\ell)$.
We assume that the group chain $\G^x_\cU$ is in the normal form, so by Lemma \ref{lemma-normalform-transitive} $C_\ell$ acts minimally on $\widetilde U_\ell$, and also on $g_\ell^{-1}(\widetilde U_\ell)$. 
  For each $m \in C_{\ell}$ define  the conjugate element $h_{\ell}^{m} = m h_{\ell}  m^{-1}$. 
Since $h_{\ell}$ acts as the identity on $\widetilde U_{\ell}$, and   $m \cdot \widetilde U_{\ell} = \widetilde U_{\ell}$ we have  $h_{\ell}^{m} \in K_{\ell}^\G$. Thus we have a collection of maps $\{h_{\ell}^{m} \mid m \in  C_{\ell} \} \subset K_{\ell}^\G \subset \G_\ell$.

Recall that we have a clopen $W' \subset g_\ell^{-1}(\widetilde U_\ell)$ such that $h_\ell|_{W'} = \id$. Since $C_\ell$ acts minimally on $g_\ell^{-1}(\widetilde U_\ell)$, there exists $m \in C_\ell$ such that $m^{-1} \cdot y_\ell \in W'$. Then 
  $$h^m_\ell( y_\ell) = m h_\ell (m^{-1} \cdot y_\ell),$$ 
  and $h^m_\ell$ fixes every points in an open neighborhood of $y_\ell$. 
  
 Now suppose the inclusion $\Upsilon^{x,\G}_c(\Phi) \to \Upsilon^{x,\G}_s(\Phi)$ is an isomorphism, then there exists an increasing  subsequence $\{ k_\ell \mid \ell \geq 1\}$ such that the inclusion $Z_{k_\ell}^\G \to K_{k_\ell}^\G$ is an isomorphism for all $\ell \geq 1$. For simplicity, set $k_\ell = \ell$. Then $h_\ell \in Z_\ell^\G$. Since $C_\ell \subset \G_\ell$, this implies that $h_\ell^{-1} m^{-1} h_\ell = m^{-1}$, and so $h^m_\ell = m h_\ell m^{-1} = h_\ell$, and so $h_\ell^m$ has a non-trivial germ at $y_\ell$. This contradicts the assertion in the previous paragraph, that $h^m_\ell$ fixes an open neighborhood of $y_\ell$. Thus,  $Z_{\ell}^\G = K_{\ell}^\G$ is impossible.
\endproof

Theorem \ref{thm-incl-countable} allows us to formulate the following criterion, which is an obstruction to the existence of a non-Hausdorff element in $\cG(\fX,\G,\Phi)$.

\begin{criterion}\label{crit-upsilon-groups}
Let $(\fX,\G,\Phi)$ be a minimal equicontinuous action.  If the inclusion of direct limit countable centralizer and stabilizer groups $\Upsilon^{\G}_c(\Phi) \to \Upsilon^\G_s(\Phi)$ is an isomorphism, then the germinal groupoid $\cG(\fX,\G,\Phi)$ is Hausdorff.
\end{criterion}

\begin{ex}
{\rm
An uncountable family of pairwise non-isomorphic wild actions of a finite index torsion-free subgroup of ${\rm SL}(n,\mZ)$, for $n \geq 3$, was constructed in \cite[Theorem 1.10]{HL2019}. It was shown in the proof of the theorem that the inclusion of profinite groups $\Upsilon_c(\Phi) \to \Upsilon_s(\Phi)$ is an isomorphism. By Corollary \ref{cor-incl-prof-count} this implies that the inclusion $\Upsilon^{\G}_c(\Phi) \to \Upsilon^{\G}_s(\Phi)$ of the direct limit countable groups is an isomorphism, and so by Theorem \ref{thm-incl-countable} these actions have Hausdorff germinal groupoids.
}
\end{ex}

\subsection{Proof of Theorem \ref{thm-main1}} The proof of item (1), which states that the property that $(\fX,\G,\Phi)$ is countably wild of finite or flat type, or that it is countably dynamically wild is an invariant of the conjugacy class of the action, is given in Theorem \ref{thm-isomupsilon-countable}. 

In item (2) if $(\fX,\G,\Phi)$ is countably wild flat type, then the map of direct limit groups $\Upsilon_c^\G(\Phi) \to \Upsilon_s^\G(\Phi)$ is an isomorphism, and by Theorem \ref{thm-incl-countable} $\G$ does not have a non-Hausdorff element. Then by Criterion \ref{crit-nonH-element} the germinal groupoid $\cG(\fX,\G,\Phi)$ must have Hausdorff topology.

Item (3) is proved in Theorem \ref{thm-incl-countable}. Items (4) and (5) are proved in Corollary \ref{cor-incl-prof-count}.

\section{Actions on rooted trees}\label{sec-trees}

In this section we recall the background on group actions on rooted trees, which is necessary for the rest of the paper. The main references here are \cite{Nekrashevych2005} or \cite{Grigorchuk2011}. We restrict to $d$-ary rooted trees as these are the ones admitting self-similar actions. In a similar way, one can define actions on more general spherically homogeneous trees, and we refer the interested reader, for instance, to \cite{GL2019} for a description and some examples.

\subsection{Equicontinuous actions on $d$-ary trees}\label{sec-tree-model} Let $d \geq 2$ be an integer. 

A $d$\emph{-ary tree} $T$ is an infinite graph without cycles, consisting of a set of vertices $V = \bigsqcup_{\ell \geq 0} V_\ell$, and a set of edges $E$ chosen as follows. In $V$, the finite set $V_\ell$, called the \emph{vertex set at level} $\ell \geq 0$, is a set of cardinality $d^\ell$. The set $V_0$ is a singleton, called the \emph{root} of the tree. Edges in $E$ join vertices in $V_{\ell+1}$ and $V_\ell$ so that every vertex in $V_{\ell+1}$ is joined by an edge to a single vertex in $V_\ell$, and every vertex in $V_\ell$ is joined by edges to $d$ vertices in $V_{\ell+1}$. If $d = 2$, then $T$ is called a \emph{binary} tree. 

An infinite path in $T$ is a sequence of vertices $\{v_\ell \mid \ell \geq 0\} \subset \prod_{\ell \geq 0} V_\ell$ such that $v_{\ell+1}$ and $v_\ell$ are joined by an edge, for $\ell \geq 0$. The boundary $\partial T$ of $T$ is the collection of all infinite paths in $T$. The space $\partial T$ is a Cantor set with the relative topology from the product topology on $ \prod_{\ell \geq 0} V_\ell $. 

We label vertices in $V$ by finite words in the alphabet $\cA = \{0,1,\ldots,d-1\}$ as follows. The root $v_0 \in V_0$ is labelled by an empty word. Vertices in $V_1$ are labelled by $\cA$. Inductively, suppose $v_\ell \in V_\ell$ is labelled by a word $w_1 \cdots w_\ell$ of length $\ell$. There are $d$ vertices $v_{\ell+1,k}$, $0 \leq k \leq d-1$ in $V_{\ell+1}$ which are joined by edges to $v_\ell$, and we label $v_{\ell+1,k}$ by the word $w k = w_1 \cdots w_\ell k$ of length $\ell+1$.  

The labelling of vertices in $V$ induces the labelling of infinite paths in $\partial T$ by infinite sequences with entries in the alphabet $\cA$. Namely, an infinite sequence $w = w_1 w_2 \cdots w_\ell \cdots$ corresponds to an infinite path in $\partial T$ which passes through the vertex in $V_\ell$ labelled by the word $w_1 \cdots w_\ell$, for $\ell \geq 1$.

\begin{defn}\label{defn-automorphism}
An automorphism $g$ of a rooted $d$-ary tree $T$ is a map of $T$ which restricts to bijections $g:V \to V$ and $g: E \to E$ with the following properties:
\begin{enumerate}
\item On each level set $V_\ell$, $\ell \geq 0$, $g$ restricts to a permutation of $V_\ell$.
\item The permutations of $V_\ell$ are compatible with the tree structure. Namely, $v_\ell \in V_\ell$ and $v_{\ell+1} \in V_{\ell+1}$ are joined by an edge if and only if $g(v_\ell) \in V_\ell$ and $g(v_{\ell+1}) \in V_{\ell+1}$ are joined by an edge; and $v_{\ell+1,1}$ and $v_{\ell+1,2}$ are joined to the same vertex in $V_{\ell}$ if and only if $g(v_{\ell+1,1})$ and $g(v_{\ell+1,2})$ are joined to the same vertex in $V_\ell$.
\end{enumerate}
\end{defn}

The group of automorphisms of a rooted $d$-ary tree is denoted by $\Aut(T)$. It is well-known that $\Aut(T) = {\rm Sym}(d) \ltimes {\rm Sym}(d) \ltimes \cdots$, the infinite wreath product of the symmetric group on $d$ elements, see \cite{BOERT1996} or \cite[Section 4.1]{Lukina2021}. 

Since the action of every $g \in \Aut(T)$ restricts on each $V_\ell$ to a permutation, which is a bijective map, the action of $g$ on $T$ induces a homeomorphism $\Phi(g): \partial T \to \partial T$ of the path space. Given two infinite words $w = w_1 w_2 \cdots$ and $w' = w_1' w_2' \cdots$, define a metric on the space $\partial T$ by
 $$D\left( w,w'\right) = {2}^{-m}, \quad \textrm{where} \quad m = \min \{\ell -1 \mid w_1 w_2 \cdots w_\ell \ne w_1' w_2' \cdots w_\ell', \, \ell \geq 1\}.$$
The action of $\Aut(T)$ on $\partial T$ is by isometries relative to the metric $D$. Thus the action $(\partial T, \G,\Phi)$, where $\G \subset \Aut(T)$ is a subgroup, is equicontinuous. 

\begin{remark}
{\rm
Let ${\bf n} = (n_1,n_2,\ldots)$ be a sequence of positive integers. A \emph{spherically homogeneous tree} $T_{\bf n}$ with \emph{spherical index} ${\bf n}$ is an infinite graph without cycles, consisting of a set of vertices $V = \bigsqcup_{\ell \geq 0} V_\ell$, and a set of edges $E$, such that $V_0$ is a singleton, and for $\ell \geq 1$ every vertex in $V_{\ell-1}$ is connected by edges to precisely $n_\ell$ vertices in $V_\ell$. Thus $V_\ell$ has cardinality $n_1 \cdots n_\ell$. Setting $n_\ell = d$ for all $\ell \geq 1$ we obtain the $d$-ary, or regular, tree, described above. Similarly to the above, the group ${\rm Aut}(T_{\bf n})$ acts on the space of infinite paths in $T_{\bf n}$ equicontinuously. Morever, it is well-known that any minimal equicontinuous action on a Cantor set is conjugate to an action of a subgroup $\G \subset {\rm Aut}(T_{\bf n})$, for a suitably chosen spherical index ${\bf n}$. The procedure of constructing $T_{\bf n}$ is explained in detail, for instance, in \cite[Section 3]{GL2019}. The choice of the spherical index ${\bf n}$ is not unique, however, the prime divisors of the elements in the sequence $\{n_1,n_2,\ldots\}$ impose restrictions on possible choices of the spherical index, and on the dynamics of the action. For instance, \cite{HL2023} studies the dynamical properties of minimal equicontinuous actions of nilpotent groups depending on the finiteness of the set of prime divisors of the elements in $\{n_1,n_2,\ldots\}$. When $n_\ell = d$ for all $\ell \geq 1$, and the action of $\G$ is \emph{self-similar}, then it is possible to describe the action of the automorphisms in $\G$ by a recursive formula, which provides a powerful tool for the study of such actions. This is explained in more detail in Section \ref{sec-selfsimilar} below.

}
\end{remark}

\subsection{Self-similarity}\label{sec-selfsimilar} The automorphism group $\Aut(T)$ of a $d$-ary tree and some of its subgroups have an interesting property called self-similarity which we describe now. We first introduce some useful notation.

\begin{defn}\label{defn-subtree}
For $d \geq 2$, let $T$ be a $d$-ary tree as in Section \ref{sec-tree-model}.
Let $v_m \in V_m$ be a vertex, and let $\{v_\ell \mid 0 \leq \ell \leq m\}$ be a finite path in $T$ from the root $v_0$ to $v_m$. A \emph{subtree} $v_mT$ of $T$ consists of the vertex set 
$$ V(v_mT) = \{v_0,v_1,\ldots, v_m\}\bigcup \left(\bigcup_{\ell \geq m+1} V_\ell \right),$$ and of all edges in $E$ which join the vertices in $V(v_mT)$. 
\end{defn}

The path space of $v_mT$ is a clopen subset of $\partial T$, consisting of paths through the vertices in $V(v_mT)$, and denoted by $\partial (v_mT)$. If $w = w_1 \cdots w_m$ is a finite word which labels $v_m$, then we may also denote the subtree by $wT$, and its path space by $\partial (w T)$. The path space $\partial (v_mT) = \partial (w  T)$ consists of all paths which start with the finite subword $w_1 \ldots w_m$. 

 Every vertex of $wT \cap V_\ell$ for $\ell \geq m$ has a label of the form $wk$, where $k$ is a word of length $\ell-m$ in the alphabet $\cA = \{0,\ldots, d-1\}$. Every letter in $w$  or $k$ is a symbol in $\cA$, so there is a bijection on the sets of vertices
  \begin{align}\label{eq-treehomeo} \pi_{w}: wT\cap V  \to V: wk \mapsto k, \end{align}
which induces a homeomorphism of path spaces $\overline{\pi_w}: \partial (wT) \to \partial T$. In the arguments below we use the word notation for the vertices in $V$.
  
Now let $g \in \Aut(T)$, and suppose $g$ maps $w \in V_m$ to a vertex $g(w) \in V_m$. 
By property (2) in Definition \ref{defn-automorphism}, the action of $g$ maps the clopen set $ \partial (wT) $ homeomorphically onto the clopen set $ \partial (g(w)T) $.
More precisely, for each vertex $wk \in wT$ there is a unique vertex $g(wk) \in g(w)T$, which is labelled by the word $g(w)k'$ for some finite word $k'$. Composing the bijections \eqref{eq-treehomeo} for $w$ and $g(w)$, we can define the bijection, called the \emph{section} of $g$ at $w$
  \begin{align}\label{eq-restriction-g}g|_w = \pi_{g(w)} \circ g \circ \pi_w^{-1}: V \to V: k \mapsto k', \end{align}
which defines an automorphism of the tree $T$, and so induces a homeomorphism $\Phi(g|_w):\partial T \to \partial T$.

 \begin{defn}\label{defn-selfsim}\cite[Definition 1.5.3]{Nekrashevych2005}
 Let $\G$ be a subgroup of $\Aut(T)$. Then $\G$ is \emph{self-similar} if for every $g \in \G$ and every vertex $w$ in $T$ the section $g|_w$ defined by \eqref{eq-restriction-g} is in $\G$.
 \end{defn}

For instance, $\Aut(T)$ itself is self-similar. 
Self-similar subgroups of $\Aut(T)$ admit the following description, which make their study especially accessible.

In the compositions of maps below, we compose on the left.

Suppose $\G \subset \Aut(T)$ is self-similar, and let $g \in \G$. Recall that $V_1$ is a set with $d$ vertices. Set $\sigma_g = g|V_1$, that is, $\sigma_g$ is a permutation of vertices in $V_1$ induced by the action of $g$. For every $w \in V_1 \cong \cA$ we have $g|_w \in \G$, so we can define a function $f_g: V_1 \to \G^{|V_1|}: w \mapsto g|_{\sigma_g^{-1}(w)}$. Then $g$ acts on $T$ as an element $(f_g, \sigma_g)$ of the semi-direct product $\G^{|V_1|} \rtimes {\rm Sym}_d$. More precisely, for $w = w_1 w_2 \cdots \in \partial T$, the element $g$ acts on $w_1$ as $\sigma_g$, and on the infinite sequence $w_2 w_3 \cdots$ as $f_g(w_1) = g|_{w_1}$. Thus we can represent $g$ as a composition
   \begin{align}\label{eq-confusingnotation}g = (g|_{\sigma^{-1}_g(0)}, g|_{\sigma^{-1}_g(1)}, \ldots, g|_{\sigma^{-1}_g(d-1)}) \circ \sigma_g, \end{align}
 where $\sigma_g = (1, \sigma_g) \in \G^{|V_1|} \rtimes {\rm Sym}_d$ and $(g|_{\sigma^{-1}_g(0)}, g|_{\sigma^{-1}_g(1)}, \ldots, g|_{\sigma^{-1}_g(d-1)}) \in \G^{|V_1|}$. Here $1$ denotes the trivial function in $\G^{|V_1|}$ which assigns to each $w \in V_1$ the identity map $1 \in \Aut(T)$. Computing \eqref{eq-confusingnotation} we first apply the permutation $\sigma_g$ to $V_1$, and then the maps $g|_{\sigma^{-1}_g(w)}$ to the subtrees $wT$, $0 \leq w \leq d-1$. 
 
 Alternatively, we can also write $g$ as the following composition (again, we compose maps on the left)
    \begin{align}\label{eq-moreconfusingnotation}g = \sigma_g \circ  (g|_{0}, g|_{1}, \ldots, g|_{d-1}), \end{align}
 that is, when computing the action of $g$ we first apply the maps $g_w$ to the subtrees $wT$, $0 \leq w \leq d-1$, and then we apply the permutation $\sigma_g$ of $V_1$. Different sources in the literature use one or the other of these two ways to write an automorphism $g \in \Aut(T)$ as a composition of two maps.  We will use \eqref{eq-confusingnotation}. Formulas \eqref{eq-confusingnotation} and \eqref{eq-moreconfusingnotation} together give the relation
  \begin{align}\label{eq-switch}  \sigma_g \circ  (g|_{0}, g|_{1}, \ldots, g|_{d-1}) = (g|_{\sigma^{-1}_g(0)}, g|_{\sigma^{-1}_g(1)}, \ldots, g|_{\sigma^{-1}_g(d-1)}) \circ \sigma_g, \end{align}
 which one can use to change from one notation to another one. 
 
 \begin{remark}\label{remark-diffnotation}
 {\rm
 Among the articles whose results we use in this paper, \cite{Pink2013,Noce2021} uses the same convention \eqref{eq-confusingnotation} for writing an element of $\Aut(T)$ as us. Paper \cite{BN2008} only considers binary trees and uses a slightly different convention: when they write $g = (g_0,g_1) \sigma$, the action of $g$ on an infinite sequence $w = w_1 w_2 \cdots$ is given by $g(0 w_2 \cdots) = 1 g_0(w_2 \cdots)$ and $g(1 w_2 \cdots) = 0 g_1(w_2 \cdots)$. In our notation, this is equivalent to \eqref{eq-moreconfusingnotation} with composition of maps on the right, and only makes a difference in the definition of the element $a_1$ in Section \ref{sec-kv}.
 }
 \end{remark}

\subsection{Non-Hausdorff elements in contracting groups}\label{sec-contracting-nonH} The self-similar property of certain subgroups of $\Aut(T)$, where $T$ is a $d$-ary tree, allows us to determine when the germinal groupoids associated to their actions have the non-Hausdorff property.
 
\begin{defn}\label{defn-contracting}
Let $\G \subset \Aut(T)$ be a self-similar group. We say that $\G$ is \emph{contracting}, if there exists a finite set $\cN \subset \G$ such that for every $g \in \G$ there is $n_g \geq 0$ such that for all finite words $w$ of length at least $n_g$ we have $g|_w \in \cN$. 
\end{defn} 

The set $\cN$ is called the \emph{nucleus} of the group $\G$, if $\cN$ is the smallest possible set satisfying Definition \ref{defn-contracting}. Iterated monodromy groups of post-critically finite polynomials in Theorem \ref{thm-polyn} and, more generally, iterated monodromy groups of post-critically finite rational functions are known to be contracting \cite[Theorem 6.4.4]{Nekrashevych2005}.

We also consider the following special subsets of the group $\G$, introduced in \cite{Jones2015}. Let
 \begin{align}\label{eq-N0} \cN_0 = \{g \in \G \mid g|_v = g \textrm{ for some non-empty word}\, w \in V \}. \end{align}
The set $\cN_0$ is always non-empty, as it contains the identity of $\G$.  It is proved in \cite[Proposition 3.5]{Jones2015} that if $\G$ is contracting, then $\cN_0$ is finite and the nucleus of $\G$ is given by
  $$\cN = \{h \in \G \mid  h = g|_w \textrm{ for some }g \in \cN_0 \textrm{ and }w \in V \}.$$
Also define
  \begin{align}\label{eq-N1} \cN_1 = \{g \in \G \mid g|_w = g \textrm{ and }g(w) = w \textrm{ for a non-empty  word } w \in V\}. \end{align} 
Then $\cN_1$ contains elements in $\cN_0$ which fix at least one path in $T$, so $\cN_1 \subset \cN_0$ and $\cN_1$ is finite for contracting actions. Also, $\cN_1$ is non-empty as it contains the identity of $\G$. The following statement is proved in the last paragraph of \cite[Section 4]{Jones2015} on p. 2033.

\begin{lemma}\cite{Jones2015}\label{lemma-N1torsion}
Let $\G \subset \Aut(T)$ be a contracting self-similar group. Then every $g \in \cN_1$ is torsion.
\end{lemma}

We will need the following statement, proved in the author's work \cite{Lukina2021}.

\begin{lemma}\label{non-hausdorff-n1}\cite[Lemma 4.5]{Lukina2021}
Let $\G \subset \Aut(T)$ be contracting, and suppose $\G$ contains a non-Hausdorff element $h$. Then there is a non-Hausdorff element $g$ in $\mathcal{N}_1$.
\end{lemma}

Summarizing the results of Lemmas \ref{lemma-N1torsion} and \ref{non-hausdorff-n1} we obtain the following property of non-Hausdorff elements in contracting self-similar groups.

\begin{criterion}\label{crit-contracting}
Let $T$ be a $d$-ary rooted tree.
Let $\G \subset \Aut(T)$ be a contracting self-similar group. If the associated germinal groupoid $\cG(\fX,\G,\Phi)$ is non-Hausdorff, then there exists a non-Hausdorff element $g \in \G$ which has finite order.
\end{criterion}

Criterion \ref{crit-contracting} is one of the main components of the proof of Theorem \ref{thm-polyn}.

\section{Germinal groupoids associated to quadratic PCF polynomials}\label{sec-PCF}

In this section we consider iterated monodromy groups associated to quadratic polynomials. We refer to \cite{Nekrashevych2005} for details on how, given a polynomial $f: \mC \to \mC$ of degree $d \geq 2$, one can define an iterated monodromy group. In our proofs below we use presentations for generators of such groups available in the literature. The goal of this section is to prove Theorems \ref{thm-no-crit} and \ref{thm-polyn}.

Groups $\fK(v)$ and $\fK(w,v)$ defined in Sections \ref{sec-kv} and \ref{sec-kwv} were introduced in \cite{BN2008}. It was proved in \cite{BN2008} that iterated monodromy groups associated to quadratic polynomials correspond to a proper subset of the set of these groups. In this section $T$ is a binary tree, and $w$ and $v$ are finite words in the alphabet $\cA = \{0,1\}$. In $\fK(v)$, the word $v$ may be empty.

\subsection{Hausdorff groupoids for groups $\fK(v)$} \label{sec-kv} Let $T$ be a binary tree. We first show that the germinal groupoid $\cG(\partial T,\fK(v),\Phi)$ associated to the action of the group $\fK(v)$ on the boundary of $T$ has Hausdorff topology.

Let $v = x_1 \cdots x_{n-1}$ be a non-empty word with $x_i \in \{0,1\}$. 

We denote by $\sigma$ an element of $\Aut(T)$ such that for each infinite sequence $z_1 z_2 \cdots \in \partial T$, $\sigma(z_1) = z_1 +1 \mod 2$, and $\sigma(z_\ell) = z_\ell$ for $\ell \geq 2$. 
As in \cite{BN2008}, define the subgroup $\fK(v) = \langle  a_1,\ldots, a_n\rangle$ of $\Aut(T)$ by
    \begin{align}\label{gen-a-v} a_1 = (a_n,1) \sigma, & & a_{i+1} = \left\{ \begin{array}{ll} (a_i,1), & \textrm{ if } x_i = 0 \\ (1, a_i), & \textrm{ if } x_i = 1\end{array} \right., \textrm{ for } 1 \leq i \leq n-1.\end{align}
 Groups $\widetilde G_r$ in item (4) of \cite[Theorem 1.5]{Lukina2021}, with $r \geq 2$, are the groups $\fK(v)$ with $x_i = 0$ for $1 \leq j \leq r-1$. The formula for $a_1$ in \eqref{gen-a-v} is different from that in \cite{BN2008} due to the slight difference in our notation, see Remark \ref{remark-diffnotation}.
 
 If $v$ is an empty word, then $\fK()$ is generated by $a_1 = (a_1,1)\sigma$, and its action is an odometer action, see for instance \cite[Example 4.2]{Lukina2021}. The odometer action on $\partial T$ is free, and so the associated germinal groupoid $\cG(\partial T, \fK(),\Phi)$ has Hausdorff topology by Criterion \ref{crit-hausdorff2}.
 
 \begin{thm}\label{thm-hausdorff-periodic}
Let $T$ be a binary tree, let $n \geq 2$ and let $v$ be a non-empty word of length $n-1$ with letters in $\{0,1\}$. Let $\fK(v)$ be a group with $n$ generators, defined by \eqref{gen-a-v}. Then the germinal groupoid $\cG(\partial T,\fK(v),\Phi)$ associated to the action of $\fK(v)$ on $\partial T$  has Hausdorff topology.
\end{thm}

\proof Suppose the germinal groupoid $\cG(\partial T, \fK(v),\Phi)$ associated to the action of $\fK(v)$ on $\partial T$ is non-Hausdorff. By \cite[Lemma 3.2]{BN2008} the group $\fK(v)$ is self-similar and contracting. Then by Criterion \ref{crit-contracting} the group $\fK(v)$ must contain a non-Hausdorff element which is torsion. But by \cite[Proposition 3.11]{BN2008} the group $\fK(v)$ is torsion-free. Thus $\cG(\partial T, \fK(v),\Phi)$ must have Hausdorff topology.
\endproof

\subsection{A countably dynamically wild action with Hausdorff germinal groupoid}

In this section we prove Theorem \ref{thm-no-crit}, giving an example of a countably dynamically wild action with Hausdorff germinal groupoid. This shows that Criterion \ref{crit-upsilon-groups} is not an `if and only if' statement. That is, if Criterion \ref{crit-upsilon-groups} is not satisfied, the germinal groupoid associated to a minimal equicontinuous action may still have Hausdorff topology.

 \begin{thm}\label{thm-periodic-Upsilon}
Let $T$ be a binary tree, and consider the self-similar contracting group
  $$\fK(1) = \langle a_1, a_2\rangle, \, \textrm{where} \, a_1 = (a_2,1) \sigma, \, a_2 = (1, a_1).$$
Then the action of $\fK(1)$ does not satisfy Criterion \ref{crit-upsilon-groups} while the germinal groupoid $\cG(\partial T, \fK(1),\Phi)$ has Hausdorff topology.
\end{thm}

\begin{figure}
\includegraphics[width=7cm]{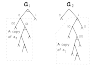}
\caption{Portaits of generators of $\fK(1)$:  we have sections $a_1|_0 = a_2$, $a_2|_1 = a_1$. Two short edges at a vertex $u$ correspond to the trivial section $a_i|_u = 1$, $i \in \{1,2\}$.}
\label{fig-s0r2-11}
\end{figure}

Portraits of the generators of $\fK(1)$ are shown in Figure \ref{fig-s0r2-11}. An arc, joining two edges emanating from the vertex labelled by the word $u = u_1 u_2 \cdots u_\ell$ means that the action of $a_i$, $i \in \{1,2\}$ changes the letter following $u_\ell$ in an infinite sequence $z = u_1u_2 \cdots u_\ell z_{\ell+1} \cdots$, i.e. $a_i(z_{\ell+1}) = z_{\ell+1} + 1\mod 2$. The action of $a_i$ may or may not change the letters before or after $z_{\ell+1}$, this depends on $a_i$ and $z$. Two short edges attached to $u$ mean that the action does not change any  letter after $u_\ell$, i.e. the corresponding section is trivial, $a_i|_u = 1 \in \Aut(T)$.
 
\proof Let $z = 1^\infty$ be an infinite path in $\partial T$ with a neighborhood basis $\cU = \{U_\ell\}_{\ell \geq 1}$, where $U_\ell = \partial (1^\ell T) $, and the subtree $1^\ell T$ is as in Definition \ref{defn-subtree}. Let $\G_\ell$ be the isotropy group of $U_\ell$, which in this case is equal to the stabilizer of the vertex $1^\ell$ in $T$.  
Let $K_\ell^\G$, $\ell \geq 0$,  be the stabilizer subgroups of $\G$ defined by \eqref{eq-KZG}, that is, $g \in K_\ell^\G$ if and only if $g|U_\ell = \id$, which is equivalent to the conditions $g(1^\ell) = 1^\ell$ and $g|_{1^\ell} = 1$ being satisfied at the same time. Let $Z_\ell^\G$ be the centralizer subgroups of $\G$ defined by \eqref{eq-KZGZ}, that is, $g \in Z_\ell^\G$ if and only if $[g,h] = 1$ for all $h \in \G_\ell$, where $1$ denotes the identity in $\Aut(T)$. To show that Criterion \ref{crit-upsilon-groups} is not satisfied we must show that the inclusion of the direct limit groups 
  $$\Upsilon_c^{z,\G} =  \varinjlim \cG(Z_{\ell}^\G, \iota^{\ell'}_{\ell} , \mN) \to  \Upsilon_s^{z,\G} = \varinjlim \cG(K_{\ell}^\G, \iota^{\ell'}_{\ell} , \mN)$$
  is proper. To this end we will show that for each $\ell > 4$, the inclusion $Z_\ell^\G \to K_\ell^\G$ is proper, that is, there exists $g_\ell \in K_\ell^\G$ and $h_\ell \in \G_\ell$ such that $g_\ell^{-1} h_\ell g_\ell \ne h_\ell$.

Denote by $1_\ell$ the trivial permutation of $V_\ell$. In the computations below, we repeatedly use formulas \eqref{eq-confusingnotation} and \eqref{eq-moreconfusingnotation}, i.e. the equality $\sigma (g_0,g_1) = (g_1,g_0) \sigma$, and the fact that $\sigma^2 = 1$, the identity in $\Aut(T)$. Also, suppose we can write for $g \in \Aut(T)$ 
  $$g = (g|_{0^\ell}, g|_{0^{\ell-1}1}\ldots, g|_{1^{\ell-1}0}, g|_{1^\ell}) 1_\ell,$$
where $g|_u$ is the section of $g$ at $u$, and $u$ is the label of a vertex in $V_\ell$, i.e. the word of length $\ell$ in the alphabet $\{0,1\}$. Then, for the square of $g$ we have
\begin{align*}g^2 = ((g|_{0^\ell})^2, (g|_{0^{\ell-1}1})^2\ldots, (g_{1^{\ell-1}0})^2, (g_{1^\ell})^2) 1_\ell,\end{align*}
For any $u$ we have either $(g|_u)|V_1 = 1_1$, or $(g|_u)|V_1 = \sigma$. Since $\sigma$ has order two, this implies that $(g|_u)^2|V_1 = 1_1$, for any $u \in V_\ell$, and so $g^2|V_{\ell+1} = 1_{\ell+1}$. 
 
So we have 
\begin{align} \label{eq-producttwo} a_1 a_2 & = (a_2,1)\sigma(1,a_1) = (a_2a_1,1)\sigma, \\ \label{eq-prodtwoinv}a_2a_1 & = (1,a_1)(a_2,1) \sigma = (a_2,a_1)\sigma.\end{align} 
Next, note that
  \begin{align*} (a_1a_2)^2 & = (a_2a_1,1)\sigma(a_2a_1,1)\sigma = (a_2a_1,a_2a_1)1_1, \\ 
       (a_2a_1)^2 & = (a_2,a_1)\sigma (a_2,a_1)\sigma = (a_2a_1,a_1a_2)1_1,\end{align*}
so
 $$(a_1a_2)^{2^2}=(a_1a_2)^4 = \left((a_2a_1)^2,(a_2a_1)^2 \right) = (a_2a_1,a_1a_2,a_2a_1,a_1a_2)1_2.$$
 Inductively, computing $a^{2 \cdot 2^{\ell-1}}$, for $\ell \geq 2$ we obtain that 
   \begin{align} \label{eq-a1a22i}(a_1a_2)^{2^\ell} = (g_{0^\ell}, \ldots, , g_{1^\ell}) 1_\ell, & & ~ \textrm{where}~ &  \left\{ \begin{array}{l}g_{1^\ell}  = a_2a_1  \, \textrm{~for~} \ell \textrm{~odd~}, \\ 
    g_{1^\ell} = a_1a_2  \, \textrm{~for~} \ell \textrm{~even~}. \end{array} \right.\end{align}

Next, $a_1^2 = (a_2,a_2)$, and $a_2^2 = (1,a_1^2) = (1,1,a_2,a_2)1_2$. We have, inductively, the sections
  $$a_2|_1 = a_1, \, ~ \,a_2^2|_{11} = a_2,$$
 and, for $\ell \geq 2$
   \begin{align}\label{eq-a22l}a_2^{2^\ell} =  (g_{0^{2\ell}}, \ldots, g_{1^{2\ell}}) 1_{2\ell}, & ~ ~\textrm{where} ~  \, 
   g_{1^{2\ell}}  = a_2.  \end{align}
Since $(a_1a_2)^{2^\ell}|V_\ell = 1_\ell$, from \eqref{eq-a1a22i} we have  the inverse    
   \begin{align*}(a_1a_2)^{-2^\ell} = (g_{0^\ell}^{-1}, g_{0^{\ell-1}1}^{-1}\ldots, g_{1^{\ell-1}0}^{-1}, g_{1^\ell}^{-1}) 1_\ell, \end{align*}
where   $(a_1a_2)^{-1} = (1,(a_2a_1)^{-1}) \sigma$, and $(a_2a_1)^{-1} = (a_1^{-1},a_2^{-1})\sigma$.
   
Let $\ell > 4$ be odd, so $\ell-1$ is even. Then $m_\ell = \frac{\ell-1}{2}$ is an integer, and $a_2^{2^{m_\ell}}|_{1^{\ell-1}} = a_2$.

Let $g_\ell = (a_1a_2)^{-2^{\ell-1}}\, a_2^{2^{m_\ell}} \, (a_1a_2)^{2^{\ell-1}}$. We compute the section of $g_\ell$ at $1^{\ell-1}$:
 \begin{align}\label{eq-conja22i} g_\ell|_{1^{\ell-1}} = (a_1a_2)^{-2^{\ell-1}}\, a_2^{2^{m_\ell}} \, (a_1a_2)^{2^{\ell-1}}|_{1^{\ell-1}} & = (a_1a_2)^{-1}\, a_2 \, (a_1a_2) \\ \nonumber & =  (1,(a_2a_1)^{-1}) \sigma (1,a_1)(a_2a_1,1)\sigma = (a_1,1),\end{align}
which implies that $g_\ell(1^\ell) = 1^\ell$ and $g_\ell|_{1^{\ell}} = 1$, and so $g_{\ell} \in K_{\ell}^\G$.

Also, since by \eqref{eq-a1a22i} $(a_1a_2)^{2^\ell}|V_{\ell} = 1_{\ell}$, its action fixes the vertex $1^\ell $. Choose $h_{\ell} = (a_1a_2)^{-2^\ell} \in \G_{\ell}$.  Note that since $(a_1a_2)^{2^{\ell-1}}(1^{\ell-1}) = 1^{\ell-1}$ and  $(a_1a_2)^{2^{\ell-1}}|_{1^{\ell-1}} = a_1a_2$ by \eqref{eq-a1a22i}, then $(a_1a_2)^{2^{\ell}}|_{1^{\ell-1}} = (a_1a_2)^2$, and $(a_1a_2)^{-2^{\ell}}|_{1^{\ell-1}} = (a_1a_2)^{-2}$.

Suppose that $g_{\ell} \in Z^\G_{\ell}$, then $[g_{\ell},h_{\ell}] = 1 \in \Aut(T)$, and we must have $h_{\ell}^{-1}g_{\ell}h_{\ell} = g_{\ell}$. We compute the section at $1^{\ell-1}$, using \eqref{eq-conja22i}:
\begin{align*} h_{\ell}^{-1}g_{\ell}h_{\ell} |_{1^{\ell-1}}& = (a_1a_2)^2(a_1a_2)^{-1}a_2(a_1a_2)(a_1a_2)^{-2} = (a_1a_2) a_2 (a_1a_2)^{-1} \\ & = (a_2a_1,1) \sigma (1,a_1)(1,(a_2a_1)^{-1})\sigma = (a_2a_1a_2^{-1},1).\end{align*}
From \eqref{eq-conja22i} we have that $h_{\ell}^{-1}g_{\ell}h_{\ell} |_{1^{\ell-1}} = g_\ell|_{1^{\ell-1}}$ if and only if $a_2a_1 = a_1a_2$, which is only possible if $a_1$ is trivial, see \eqref{eq-producttwo} - \eqref{eq-prodtwoinv}. Therefore, $h_{\ell}^{-1}g_{\ell}h_{\ell} |_{1^{\ell-1}} \ne g_\ell|_{1^{\ell-1}}$, and $g_\ell \notin Z^\G_\ell$.  Thus, for $\ell \geq 4$ odd the inclusion $Z_\ell^\G \to K_\ell^\G$ is proper.

Now suppose $\ell > 4$ is even, and so $\ell-1$ is odd, and $\ell-2$ is even. Then $m_\ell = \frac{\ell-2}{2}$ is an integer, and $a_2^{2^{m_\ell}}|_{1^{\ell-2}} = a_2$. We choose $g_\ell = (a_1a_2)^{-2^{\ell-2}}\, a_2^{2^{m_\ell}} \, (a_1a_2)^{2^{\ell-2}}$, then by \eqref{eq-conja22i}
 \begin{align}\label{eq-conja22iodd} g_\ell|_{1^{\ell-2}} = (a_1a_2)^{-2^{\ell-2}}\, a_2^{2^{m_\ell}} \, (a_1a_2)^{2^{\ell-2}}|_{1^{\ell-2}} & = (a_1,1),\end{align}
 so $g_\ell|U_{\ell-1} = \id$. Since $U_\ell \subset U_{\ell-1}$, we have $g_\ell|U_\ell = \id$ and $g_\ell \in K^\G_\ell$.
 
Choose $h_\ell = (a_1a_2)^{-2^{\ell}}$, then $h_\ell \in \G_\ell$. Then, computing the section of $h_\ell^{-1} g_\ell h_\ell$ at $1^{\ell-2}$ we obtain
 \begin{align*} h_{\ell}^{-1}g_{\ell}h_{\ell} |_{1^{\ell-1}}& = (a_1a_2)^4(a_1a_2)^{-1}a_2(a_1a_2)(a_1a_2)^{-4}  = ((a_2a_1)^2a_1(a_2a_1)^{-1},1),\end{align*}
which is only equal to $g_\ell|_{1^{\ell-1}}$ if $(a_2a_1)^2a_1 = a_1(a_2a_1)^2$. The latter is true if and only if $a_1a_2 = a_2a_1$, which, in its turn, holds only if $a_1$ is trivial, which is not the case. Thus $h_\ell^{-1} g_\ell h_\ell \ne g_\ell$, and for even $\ell > 4$ the inclusion $Z_\ell^\G \to K_\ell^\G$ is proper. It follows that the inclusion of direct limit groups $\Upsilon_c^{z,\G} \to \Upsilon_s^{z,\G}$ is proper, and Criterion \ref{crit-upsilon-groups} is not satisfied. However, by Theorem \ref{thm-hausdorff-periodic} the germinal groupoid associated to the action has Hausdorff topology.
 
\endproof

\subsection{Non-Hausdorff groupoids for groups $\fK(w,v)$}\label{sec-kwv}

We show that, for the action of a group $\fK(w,v)$ on the binary tree $T$, the associated germinal groupoid in almost all cases has non-Hausdorff topology, by exhibiting a non-Hausdorff element in $\fK(w,v)$.

Let $w = y_1 \cdots y_k$ and $v = x_1 \cdots x_n$ be non-empty words with $y_j,x_i \in \{0,1\}$ and such that $y_k \ne x_n$. 

Recall that we denote by $\sigma$ an element of $\Aut(T)$ such that for each infinite sequence $z_1 z_2  \cdots \in \partial T$, $\sigma(z_1) = z_1 +1 \mod 2$, and $\sigma(z_\ell) = z_\ell$ for $\ell \geq 2$. 

Following \cite{BN2008}, define the subgroup of $\Aut(T)$
  $$\fK(w,v) = \langle  b_1,\ldots,b_k, a_1,\ldots, a_n\rangle,$$
where
  \begin{align}\label{gen-b}b_1 = \sigma, \, b_{j+1} = \left\{ \begin{array}{ll} (b_j,1), & \textrm{ if } y_j = 0 \\ (1, b_j), & \textrm{ if } y_j = 1\end{array} \right., \textrm{ for } 1 \leq j \leq k-1,\end{align}

\begin{align}\label{gen-a1}  a_{1} = \left\{ \begin{array}{ll} (b_k,a_n), & \textrm{ if } y_k = 0 \textrm{ and } x_n = 1 \\ (a_n, b_k), & \textrm{ if } y_k = 1 \textrm{ and } x_n = 0\end{array} \right., \end{align}
 and
   \begin{align}\label{gen-a}  a_{i+1} = \left\{ \begin{array}{ll} (a_i,1), & \textrm{ if } x_i = 0 \\ (1, a_i), & \textrm{ if } x_i = 1\end{array} \right., \textrm{ for } 1 \leq i \leq n-1.\end{align}
 For instance, the groups $\widetilde H_r$ in item (4) of \cite[Theorem 1.5]{Lukina2021}, with $k+n = r \geq 3$, are the groups $\fK(w,v)$ with $y_j = 0$ for $1 \leq j \leq k$, $x_{i} = 0$ for $1 \leq i \leq n-1$ and $x_{n} = 1$.

\begin{figure}
\includegraphics[width=12cm]{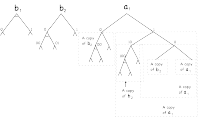}
\caption{Portaits of generators of $\fK(00,1)$:  we have sections $a_1|_0 = b_2$, $a_1|_{00} = b_1$, $a_1|_{01} = 1$, $a_1|_{1^\ell} = a_1$, $a_1|_{1^\ell0} = b_2$ for $\ell \geq 1$. Two short edges at a vertex $u$ correspond to the trivial section $g|_u = 1$, for $g \in \{b_1,b_2,a_1\}$.}
\label{fig-s0r2-111}
\end{figure}

As an example, Figure \ref{fig-s0r2-111} contains the portraits of generators of the group $\fK(00,1)$. An arc, joining two edges emanating from a vertex labelled by a finite word $u = u_1 u_2 \cdots u_\ell$ means that the action of $g$, where $g \in \{b_1,b_2,a_1\}$, changes the letter following $u_\ell$ in an infinite sequence $z = u_1u_2 \cdots u_\ell z_{\ell+1} \cdots$, i.e. $g(z_{\ell+1}) = z_{\ell+1} + 1\mod 2$. The action of $g$ may or may not change the letters before or after $z_{\ell+1}$, this depends on $g$ and $z$. Two short edges attached to $u$ mean that the action does not change any letter after $u_{\ell}$, i.e. the section at $u$ is trivial, $g|_u = 1 \in \Aut(T)$.

\begin{thm}\label{thm-thm12}
Let $T$ be a binary tree. Let $k,n \in \mN$ such that $k \geq 1$, $n \geq 1$, and $k+n \geq 3$. Then the group $\fK(w,v)$ generated by \eqref{gen-b} - \eqref{gen-a} contains a non-Hausdorff element, and so the germinal groupoid associated to the action of $\fK(w,v)$ on $\partial T$ has non-Hausdorff topology.
\end{thm}

\proof

We will consider two cases, when $n=1$ and when $n \geq 2$. In each case, we will find a non-Hausdorff element in the acting group $\fK(w,v)$. Then by Criterion \ref{crit-nonH-element} the topology of the germinal groupoid $\cG(\partial T,\fK(w,v),\Phi)$ is non-Hausdorff. 

Recall that for a finite word $u$ in the alphabet $\cA =\{0,1\}$, $uT$ denotes the subtree of $T$ containing all paths passing through the vertex $u$, see Definition \ref{defn-subtree}. These paths form a clopen subset of $\partial T$, denoted by $\partial (uT)$ and consisting of infinite sequences starting with the word $u$. For an element $g \in \Aut(T)$, the notation $g|\partial (uT) = \id$ means that $g$ acts on the clopen subset $\partial (uT)$ as the identity map. The condition $g|\partial (uT) = \id$ is equivalent to the conditions $g(u) = u$ and $g|_u = 1 \in \Aut(T)$ holding at the same time, where $g|_u$ denotes the section of $g$ at $u$ as in \eqref{eq-restriction-g}. If $g|_u = 1$ and $g(u) \ne u$, then $g|\partial (uT) \ne \id$. In this case, given an infinite sequence $uz = u_1 \cdots u_\ell z_{\ell+1} \cdots$, the action of $g$ preserves all letters $z_j$, $j \geq \ell+1$, but acts non-trivially on the letters in the word $u = u_1 \cdots u_\ell$.

\begin{lemma}\label{pre-periodic-fixedpoint}
In Theorem \ref{thm-thm12}, suppose that $n = 1$. Then $a_1 = a_n$ is non-Hausdorff.
\end{lemma}

\proof We will find an infinite path $z = z_1z_2 \cdots \in \partial T$, such that $a_n(z) = z$, a descending collection of clopen neighborhoods $\{W_\ell\}_{\ell \geq 1}$ with $\bigcap_{\ell \geq 1} W_\ell = \{z\}$ and, for each $\ell \geq 1$, a clopen subset $O_\ell \subset W_\ell$, such that the restriction $a_1|O_\ell$ is the identity map, while the restriction $a_1|W_\ell$ is non-trivial. 

For a symbol $u \in \{0,1\}$, denote by $u' = u+1 \mod 2$, i.e. if $u = 0$ then $u' = 1$, and if $u = 1$ then $u' = 0$.

If $n = 1$, then $k \geq 2$. In particular, we have $b_k = (b_{k-1},1)$ if $y_{k-1} = 0$, or $b_k = (1,b_{k-1})$ if $y_{k-1} = 1$. Then
$b_k|\partial(y_{k-1}' T) = \id$, and $b_k|\partial(y_{k-1} T)$ is non-trivial. Equivalently, $b_k|V_1 = 1_1$, that is, $b_k(y_{k-1}) = y_{k-1}$ and $b_k(y_{k-1}') = y_{k-1}'$, and, using the notation for sections, $b_k|_{y_{k-1}'} = 1$ while $b_k|_{y_{k-1}} = b_{k-1}$ is non-trivial.

Note that $a_1|V_1 = 1_1$ by definition, where $1_\ell$ denotes the trivial permutation of $V_\ell$ for $\ell \geq 0$. We set $w_1 = x_n = x_1$,  $W_1 = \partial T$, and $O_1 = \partial (x_1' y_{k-1}' T)$. By the remarks above, $a_1(x_1'y_{k-1}') = x_1' y_{k-1}'$,   and
$$a_1|_{x_1'y_{k-1}'} =  b_k|_{y_{k-1}'} = 1, $$
so $a_1|O_1 = \id$. The restriction $a_1|W_1$ is clearly non-trivial.

Similarly, for $\ell \geq 2$ let $z_\ell = x_n$, and set $W_\ell = \partial (z_1 \cdots z_{\ell-1} T)$. Then, using the recursive definition of $a_1$ we have $a_1(z_1 \cdots z_{\ell-1}) = z_1 \cdots z_{\ell-1}$, and
  $$a_1 |_{z_1 \cdots z_{\ell-1}} = a_1, \quad \quad a_1 |_{z_1 \cdots z_{\ell-1} x_n'} = b_k,$$
 so $a_1|W_\ell$ is non-trivial. Now set $O_\ell =  \partial (z_1 \cdots z_{\ell-1} x_n' y_{k-1}' T)$, then 
  $O_\ell \subset W_\ell$, and the restriction $a_1|O_\ell$ is the trivial map, since 
   $$a_1(z_1 \cdots z_{\ell-1} x_n' y_{k-1}' ) = z_1 \cdots z_{\ell-1} a_1(x_n' y_{k-1}') = z_1 \cdots z_{\ell-1} x_n' b_k(y_{k-1}') = z_1 \cdots z_{\ell-1} x_n' y_{k-1}' , $$
 and  
   $$a_1|_{z_1 \cdots z_{\ell-1} x_n' y_{k-1}' } = b_k|_{y_{k-1}' } = 1.$$
Since $a_1(z_1 \cdots z_\ell) = z_1 \cdots z_\ell$ for all $\ell \geq 1$,  $a_1$ fixes every letter $z_\ell$ in $z$, and therefore $a_1(z) = z$.   
We have shown that $a_1$, $z$, $\{W_\ell \mid \ell \geq 1\}$ and $\{O_\ell \mid \ell \geq 1\}$ as chosen above, satisfy Definition \ref{defn-nonH-rev}. Therefore, $a_1$ is a non-Hausdorff element.
\endproof

\begin{lemma}\label{pre-periodic-orbit}
In Theorem \ref{thm-thm12}, suppose that $n \geq 2$. Then for $1 \leq i \leq n$, the element $a_i$ is non-Hausdorff.
\end{lemma}

\proof First we show that $a_n$ is non-Hausdorff.  We will find an infinite path $z = z_1z_2 \cdots \in \partial T$, such that $a_n(z) = z$, a descending collection of clopen neighborhoods $\{W_\ell\}_{\ell \geq 1}$ with $\bigcap_{\ell \geq 1} W_\ell = \{z\}$ and, for each $\ell \geq 1$, a clopen subset $O_\ell \subset W_\ell$, such that the restriction $a_n|O_\ell$ is the identity map, while the restriction $a_n|W_\ell$ is non-trivial. 

Since $n \geq 2$, then $a_n = (a_{n-1},1)$ if $x_{n-1}=0$ or $a_n = (1,a_{n-1})$ if $x_{n-1} = 1$. Then $a_n|V_1 = 1_1$, where $1_\ell$ is the trivial permutation of $V_\ell$, $a_n|_{x_{n-1}'} = 1$, and $a_n|_{x_{n-1}} = a_{n-1}$ is non-trivial.

We set $z_1 = x_{n-1}$,  $W_1 = \partial T$, and $O_1 = \partial (x_{n-1}' T)$. Then  $a_1|O_1 = \id$.

Next, for $2 \leq j \leq n-1$ let $z_j = x_{n-j}$. Then $a_n(z_1 \cdots z_{n-1}) = a_n(x_{n-1}x_{n-2} \cdots x_1) = z_1 \cdots z_{n-1}$, and
  $$a_n|_{z_1 \cdots z_{n-1}} = a_n|_{x_{n-1} \cdots x_{1}} = a_1.$$
Further set $z_n = x_n$, then $a_n(z_1 \cdots z_{n}) = z_1 \cdots z_{n} $ and
  $$a_n|_{z_1 \cdots z_{n}}  = a_1|_{z_n } = a_1|_{x_n} =  a_n .$$  
Define $W_2 = \partial  (z_1 \cdots z_{n} T)$, and $O_2 = \partial  (z_1 \cdots z_{n} x_{n-1}'T)$, then $a_n|O_n = \id$, since $a_n(z_1 \cdots z_{n} x_{n-1}') = z_1 \cdots z_{n} x_{n-1}'$ and $a_n|_{z_1 \cdots z_{n} x_{n-1}'} = 1$.

For a word $z_1 \cdots z_n$ and a positive integer $\ell \geq 1$, denote by $[z_1 \cdots z_n]^\ell$ the concatenation of $\ell$ repetitions of $z_1 \cdots z_n$. Then define, for $\ell \geq 2$ and $1 \leq j \leq n - 1$
  $$z_{(\ell - 1) n + j} = x_{n - j}, \quad \textrm{and} \quad z_{\ell  n} = x_n.$$
Choose, for $\ell \geq 2$,  the clopen sets
   $$W_\ell = \partial \left( z_1 \cdots z_{n(\ell - 1)} T \right) = \partial \left( [x_{n-1} \cdots x_1 x_n]^{\ell - 1} T \right),  $$
 and 
   $$ O_\ell = \partial \left( z_1 \cdots z_{n(\ell - 1)}x_{n-1}' T \right) = \partial \left( [x_{n-1} \cdots x_1 x_n]^{\ell - 1}x_{n-1}' T \right),$$
 where $a_n([x_{n-1} \cdots x_1 x_n]^{\ell - 1} x_{n-1}') = [x_{n-1} \cdots x_1 x_n]^{\ell - 1} x_{n-1}'$,
  \begin{align*} ~a_n|_{[x_{n-1} \cdots x_1 x_n]^{\ell - 1}} = a_n, ~\textrm{and}~  a_n|_{[x_{n-1} \cdots x_1 x_n ]^{\ell - 1} x_{n-1}'} = a_n|_{x_{n-1}'} = 1.\end{align*}  
Then $W_\ell$ is a clopen neighborhood of $z$, with $\bigcap_{\ell \geq 1} W_\ell = z$, $O_\ell \subset W_\ell$, $a_n|W_\ell$ is non-trivial and $a_n|O_\ell = \id$, for $\ell \geq 1$. Since $a_n$ fixes every finite word $[ x_{n-1} \cdots  x_1 x_n]^\ell$,  $a_n$ fixes $z$. We have shown that $a_n$ is a non-Hausdorff element.

Now consider $a_i$, with $1 \leq i \leq n-1$. Similarly to above, we will find an infinite sequence $z^{(i)}$, and neighborhoods $W^{(i)}_\ell$ and $O^{(i)}_\ell$, for $\ell \geq 1$, such that $a_i(z^{(i)}) = z^{(i)}$, $\bigcap W^{(i)}_\ell = z^{(i)}$, $O^{(i)}_\ell \subset W^{(i)}_\ell$ and $a_i| O^{(i)}_\ell = \id$.

If $i \geq 2$, for $1 \leq j\leq i-1$ define $z^{(i)}_j = x_{i - j}$. Then, for any $i \geq 1$, set $z^{(i)}_i = x_n$. Then set
  $$W_1 = \partial (z^{(i)}_1 \cdots z^{(i)}_i T) = \partial (x_{i-1} x_{i-2} \cdots x_1 x_n T)$$ 
 and 
   $$O_1 = \partial (z^{(i)}_1 \cdots z^{(i)}_i  x_{n-1}' T) = \partial (x_{i-1} x_{i-2} \cdots x_1 x_n x_{n-1}'T).$$
Then $a_i|_{x_{i-1} \cdots x_1x_n} = a_n$, and so the action of $a_i$ restricted to $W_1$ is non-trivial. Also, $a_i({x_{i-1} \cdots x_1x_n x_{n-1}'}) =  {x_{i-1} \cdots x_1x_n x_{n-1}'}$ and $a_i|_{x_{i-1} \cdots x_1x_n x_{n-1}'} = a_n|_{x_{n-1}' } = 1$, so $a_i|O_1$ is the identity homeomorphism. We further define, for $\ell \geq 2$ and $1 \leq j \leq n-1$,
  $$z^{(i)}_{i + \ell(n-1)+j} = x_{n-j}, \quad \textrm{and} \quad z^{(i)}_{i + \ell n} = x_n,$$
 and for $\ell \geq 2$,  the clopen sets
   $$W^{(i)}_\ell = \partial \left( z_1 \cdots z_{i+ n(\ell - 1)} T \right) = \partial \left( x_{i-1} \cdots x_1 x_n [x_{n-1} \cdots x_1 x_n]^{\ell - 1} T \right), $$
 and 
   $$ O^{(i)}_\ell = \partial \left( [z_1 \cdots z_{i+n(\ell - 1)}x_{n-1}' T \right) = \partial \left( x_{i-1} \cdots x_1 x_n [x_{n-1} \cdots x_1 x_n]^{\ell - 1}x_{n-1}' T \right),$$ 
   then $a_n({x_{i-1}\cdots x_1 x_n[x_{n-1} \cdots x_1 x_n ]^{\ell - 1}x_{n-1}'}) = x_{i-1}\cdots x_1 x_n[x_{n-1} \cdots x_1 x_n ]^{\ell - 1}x_{n-1}'$, and
     \begin{align*} ~a_i&|_{x_{i-1}\cdots x_1 x_n [x_{n-1} \cdots x_1 x_n]^{\ell - 1}} = a_n, \\  a_n&|_{x_{i-1}\cdots x_1 x_n[x_{n-1} \cdots x_1 x_n ]^{\ell - 1}x_{n-1}'} = a_n|_{x_{n-1}'} = 1.\end{align*}  

   Then $W^{(i)}_\ell$ is a clopen neighborhood of $z^{(i)}$, with $\bigcap_{\ell \geq 1} W^{(i)}_\ell = z^{(i)}$, $O^{(i)}_\ell \subset W^{(i)}_\ell$, $a_i$ is non-trivial on $W_\ell^{(i)}$ and trivial on $O^{(i)}_\ell$, for $\ell \geq 1$. Since $a_i$ fixes every finite word $x_{i-1} \cdots x_1 x_n [ x_{n-1} \cdots  x_1 x_n]^\ell$, $a_i$ fixes $z^{(i)}$. We have shown that $a_i$ is a non-Hausdorff element.
\endproof

This finishes the proof of Theorem \ref{thm-thm12}.
\endproof

\subsection{Proof of Theorem \ref{thm-polyn}} Let $f: \mC \to \mC$ be a quadratic PCF polynomial. By \cite[Section 5.2]{BN2008}, if the post-critical set $P_c$ is finite and consists of a single periodic orbit, then its iterated monodromy group is the group $\fK(v)$ for some choice of $v$ ($v$ may be empty). If the post-critical set $P_c$ is a pre-periodic orbit with non-trivial pre-periodic part, then the iterated monodromy group of $f(x)$ is $\fK(w,v)$, for some choices of $v$ and $w$.

By Theorem \ref{thm-hausdorff-periodic} and the remark just before, the germinal groupoid $\cG(\partial T, \fK(v),\Phi)$ has Hausdorff topology for any choice of $v$. For $\fK(w,v)$, if the sum of lengths of $v$ and $w$ is at least $3$, then by Theorem \ref{thm-thm12} $\cG(\partial T, \fK(w,v),\Phi)$ has non-Hausdorff topology. 

The last case to consider is the case when $v$ is a word of length $1$ and $w$ is a word of length $1$. In this case, by \cite[Proposition 3.4.2]{Pink2013} the closure of $\fK(w,v)$ in $\Aut(T)$ is conjugate by an element of $\Aut(T)$ to the closure of the iterated monodromy group associated to the quadratic Chebyshev polynomial, and by \cite[Theorem 1.5(3)]{Lukina2021} the action of the closure $\overline{\fK(w,v)} \subset \Homeo(\partial T)$ on $\partial T$ is stable. Then the action $(\partial T,\fK(w,v),\Phi)$ is locally quasi-analytic, and by Criterion \ref{crit-hausdorff2} the germinal groupoid $\cG(\partial T, \fK(w,v),\Phi)$ has Hausdorff topology. This finishes the proof of Theorem \ref{thm-polyn}.  

\begin{remark}\label{remark-bounded-amenable}
{\rm
We note that the groups $\fK(v)$ and $\fK(w,v)$ described above are amenable. These groups are generated by bounded automata \cite{BN2008}, and they are amenable by the result in \cite{BKN2010}.
}
\end{remark}

\section{Non-contracting weakly branch groups of Noce}\label{sec-Md}

In this section we prove Theorem \ref{thm-Md-nonHausd}, namely, that the germinal groupoid $\cG(\partial T, \cM(d),\Phi)$, associated to the action of a weakly branch non-contracting group $\cM(d)$, $d \geq 2$, of automorphisms of a $d$-ary tree $T$, constructed in the paper by Noce \cite{Noce2021}, has non-Hausdorff topology, except when $d=2$.

 For $d \geq 2$, let $T$ be a $d$-ary tree. We use the wreath product notation as in \eqref{eq-confusingnotation} to define the generators of $\cM(d)$. Note that \cite{Noce2021} uses the alphabet $\cA' = \{1,\ldots,d\}$ in their notation, while we use the alphabet $\cA = \{0,1,\ldots,d-1\}$. Thus a group $\cM(d)$, for $d \geq 2$, has the generators
\begin{align}\label{eq-Md-groups}
m_1 & = (1,\ldots,1,m_1)(0 \ldots d-1), \\ \nonumber m_2 & = (1,\ldots,1,m_2,1)(0 \ldots d-2), \\ \nonumber & \cdots \\  
\nonumber m_{d-1} & = (1,m_{d-1},1,\ldots, 1)(01), \\ \nonumber m_d & = (m_1,m_2,\ldots,m_d),\end{align}
where $1 \in \Aut(T)$ is the identity.

\proof \emph{(of Theorem \ref{thm-Md-nonHausd})}. We want to prove that the groups $\cM(d)$, for $d \geq 3$, contain a non-Hausdorff element. Note that by \cite{Noce2021} $\cM(d)$, for $d \geq 2$, are non-contracting weakly branch groups, so Criterion \ref{crit-contracting} does not apply. We note that the non-Hausdorff elements in $\cM(d)$ which we exhibit below have infinite order.

Fix $d \geq 3$. We will show that the generator $m_d$ in \eqref{eq-Md-groups} is a non-Hausdorff element of infinite order. For that we will find an infinite path $z \in \partial T$, such that $m_d(z) = z$, a descending collection of clopen neighborhoods $\{W_\ell\}_{\ell \geq 1}$ with $\bigcap_{\ell \geq 1} W_\ell = \{z\}$, and, for each $\ell \geq 1$, a clopen subset $O_\ell \subset W_\ell$, such that the restriction $m_d|O_\ell$ is the identity homeomorphism, while the restriction $m_d|W_\ell$ is non-trivial. 

Set $z_1 = (d-1)$, $W_1 = \partial T$, and let $O_1 = \partial \left( 1(d-1)T \right)$. Then $z \in W_1$ and $O_1 \subset W_1$. Note that $m_d|V_1 = 1$, and that $m_2$ is trivial on $\partial((d-1)T)$. Then, for the section of $m_d$ at $1(d-1)$, we have
  $$m_d|_{1(d-1)} =  m_2|_{(d-1)} = 1,$$
so $m_d|O_1 = \id$.  

For a finite word $u$, denote by $[u]^\ell$ the concatenation of $\ell$ repetitions of $u$. For $\ell \geq 2$ set $z_\ell = (d-1)$, and $W_\ell = [(d-1)]^{\ell-1} \partial T$. Note that $m_d|_{(d-1)} = m_d$, and, inductively, $m_d([(d-1)]^{\ell-1}) = [(d-1)]^\ell$ and $m_d|_{[(d-1)]^{\ell-1}} = m_d$, so the action of $m_d$ on $W_\ell$ is non-trivial.

Now define $O_\ell = \partial \left( z_1 \cdots z_{\ell-1}1(d-1)T \right),$
and compute that
  $$m_d|_{z_1 \cdots z_{\ell-1}1(d-1)} = m_d|_{[(d-1)]^{\ell-1}1(d-1)}= m_d|_{1(d-1)} = m_2|_{(d-1)} = 1,$$
and also $m_d([(d-1)]^{\ell-1}1(d-1)) = [(d-1)]^{\ell-1}1(d-1)$.  Thus $m_d|O_\ell = \id$ for all $\ell \geq 1$.  Finally, $m_d (z) = z$ since $m_d$ fixes every finite word $(d-1)^\ell$. 
We have shown that $m_d$ is non-Hausdorff.  
\endproof

The group $\mathcal{M}(2)$ is also known as the \emph{long-range group} in literature \cite{AAV2013}. We show that the germinal groupoid associated to the action of this group has Hausdorff topology.\footnote{The author thanks the anonymous referee for suggesting the idea of the proof of Lemma \ref{lemma-M2Hausdorff}}.

\begin{lemma}\label{lemma-M2Hausdorff}
The germinal groupoid $\cG(\partial T, \cM(2),\Phi)$ is Hausdorff.
\end{lemma}

\proof By Criterion \ref{crit-nonH-element} the germinal groupoid associated to the action $(\partial T, \cM(2),\Phi)$ is non-Hausdorff if and only if $\cM(2)$ has a non-Hausdorff element. Suppose $g \in \cM(2)$ is non-Hausdorff, then, in particular, it has a fixed point $x = x_1x_2 \cdots $. The element $m_1 = (1,m_1) (01)$, is an odometer, and so the cyclic group generated by $m_1$ acts freely on $\partial T$. Thus $g$ can only be non-Hausdorff if its representation as a product of generators includes $m_2$.

Considering the element $m_2 = (m_1,m_2)$ we note that for any finite word $w \ne 1 \cdots 1$, i.e. $w$ is not a concatenation of only $1$'s, we have $m_2|_w = m_1$, and so the only fixed point of any power $m_2^k$, $k \in \mathbb{Z}$, is the path $z$ represented by the infinite sequence of $1$'s. Since for every $w \ne 1 \cdots 1$ we have $m_2|_w = m_1$, and $g$ is represented by a finite word in $m_1$ and $m_2$,  one can find the numbers $n_1,n_2 \in \mathbb{Z}$ such that $x = m_1^{-n_1} m_2^{n_2} m_1^{n_1} (z)$. Then $\widehat g = (m_1^{-n_1}m_2^{n_2} m_1^{n_1})^{-1} g (m_1^{-n_1}m_2^{n_2} m_1^{n_1})$ has a fixed point at $z$, and $g$ is non-Hausdorff at $x$ if and only if $\widehat g$ is non-Hausdorff at $z$. So without loss of generality we may assume that $g$ has a fixed point at $z$. Since $m_1$ acts freely on $\partial T$, we also have $g = h m_2^{\alpha}$, where $\alpha \in \{-1,1\}$ and $h \in \cM(2)$. Then $h$ also has a fixed point $z \in \partial T$ and, arguing inductively and using that $g$ is the product of a finite number of generators $m_1$ and $m_2$, we obtain that $g = m_2^k$ for some $k \in \mathbb{Z}$. This implies that, for any finite word $z_n = 1 \cdots 1$, a concatenation of $n$ copies of $1$, we have $g|_{z_n} = m_2^k$, and there exists an open neighborhood of $z$ in $\partial T$ that does not have any fixed points of $g$ except $z$. Therefore, $g$ is not a non-Hausdorff element. It follows that $\cG(\partial T, \cM(2),\Phi)$ has Hausdorff topology.
\endproof

\begin{thm}\label{exp-activity}
Let $d \geq 2$ be an integer. Then the group $\cM(d)$, generated by \eqref{eq-Md-groups}, is amenable.
\end{thm}

\proof We show that for each $d \geq 2$, the group $\cM(d)$ is generated by automata with linear activity, and, therefore, it is amenable by \cite[Theorem 1]{AAV2013}. 

Given a $d$-ary tree $T$ and $g\in \Aut(T)$, consider the wreath product representation of $g$ as in \eqref{eq-confusingnotation} and, for $n \geq 1$, the number of vertices $u$ at level $n$ such that the section $g|_u$ is non-trivial. If this number remains bounded by a constant as $n$ grows, then the automaton generating $g$ is called \emph{bounded} \cite[p.710]{AAV2013}. If this number is bounded by a function $Cn$, for some constant $C >0$, then the automaton generating $g$ has linear activity  \cite[p.710]{AAV2013}. In the group $\cM(d)$ for $d \geq 2$,  $m_1,\ldots,m_{d-1}$ are readily seen to be bounded, with only one non-trivial section at each level $n \geq 1$. However, for $m_d$ the number of non-trivial sections at level $n \geq 1$ is given  by $n(d-1)+1$, and so it grows linearly. The set of automorphisms of $\partial T$ with linear activity forms a group \cite[Section 3.8.2]{Nekrashevych2005}, therefore, all elements in $\cM(d)$ have at most linear activity. Then by \cite[Theorem 1]{AAV2013} $\cM(d)$ is amenable.
\endproof

\begin{remark}\label{remark-amenable-nonHausd}
{\rm
For $d \geq 3$, the group $\cM(d)$ is amenable by Theorem \ref{exp-activity}, and the germinal groupoid $\cG(\partial T, \cM(d),\Phi)$, associated to the action of $\cM(d)$ on the boundary of the $d$-ary tree, has non-Hausdorff topology. Thus non-contracting weakly branch groups of Noce \cite{Noce2021} provide a class of examples of actions of amenable groups with associated non-Hausdorff germinal groupoids.

Another example is given by groups $\fK(w,v)$ with $k+n \geq 3$, where $k$ is the length of $w$ and $n$ is the length of $n$. These groups are amenable since they are generated by bounded automata, and they have non-Hausdorff germinal groupoids by Theorem \ref{thm-thm12}.
}
\end{remark}


\end{document}